\documentclass{article}
\usepackage[english]{babel}
\usepackage{bussproofs,coqdoc_mod}
\usepackage{hyperref}
\usepackage[square,super,comma]{natbib}
\usepackage{textgreek,tipa}

\usepackage[utf8x]{inputenc}
\usepackage[T1]{fontenc}
\usepackage{amsmath,amsfonts,amsthm,amssymb,fancyhdr,enumerate}
\usepackage{geometry}
\usepackage{stmaryrd,bbm,xspace}

\DeclareUnicodeCharacter{8872}{\ensuremath\vDash}
\DeclareUnicodeCharacter{8594}{\ensuremath\rightharpoonup}

\DeclareMathOperator{\Hc}{Hc}

\newenvironment{code}%
 {\small\begin{list}{}%
         {\setlength{\leftmargin}{10mm}}%
         \item[]%
 }
 {\end{list}}
\newcommand{\cnf}{_{\text{CNF}}}
\newcommand{\nnf}{_{\text{NNF}}}

\makeatletter
\date{May 9, 2014}
\makeatother

\begin{document}
\title{Propositional Calculus in Coq}
\author{Floris van Doorn}
\maketitle

\begin{abstract}
    I formalize important theorems about classical propositional logic in the proof assistant Coq. The main theorems I prove are (1) the soundness and completeness of natural deduction calculus, (2) the equivalence between natural deduction calculus, Hilbert systems and sequent calculus and (3) cut elimination for sequent calculus.
\end{abstract}

\section{Introduction}\label{sec_intro}
Proof assistants (or interactive theorem provers) are computer programs which help to formalize and check the correctness of proofs. Proof assistants are used for two things. First, they are used for formal verification of hardware or software, which means proving that the hardware or software conforms to the specifications. Second, they are used to formally verify mathematical theorems. I will only discuss the second use in this paper.

To prove a mathematical theorem in a proof assistant, the user has to write down the proof in a language that the proof assistant can parse, and then the proof assistant will check correctness. Usually the proof assistant provides small amounts of automation, so that the user does not have to write down every detail. Still, a user typically has to give a lot more details to a proof assistant than is usually given on a paper proof. This means that a formalization of a mathematical theorem is much more time consuming than proving the theorem on paper, hence most research mathematics is not formalized. Still, there are some deep theorems which have been fully formalized in proof assistants, such as the four color theorem\citep{gonthier2005fourcolour} or Feit-Thompson's odd order theorem\citep{Gonthier2013Feit}, which were both formalised by Georges Gonthier in the proof assistant Coq. It is very rare that a new result is immediately accompanied by a formalization in a proof assistant, the author knows only two instances where this has been done.\citep{siles2012pure,vandoorn2013}

The goal of this project is to formalize some theorems in proof theory in a proof assistant. I chose to use the proof assistant Coq\citep{Coq:manual} for this. Coq is a proof assistant based on the predicative Calculus of Inductive Constructions, which uses the expressive power of dependent type theory as language. It uses constructive reasoning, which is very suitable for the meta-logic of proof theory, since most proofs can be given constrictively.

I have formalized the following three theorems.
\begin{enumerate}
\item The completeness theorem of classical natural deduction calculus.
\item The equivalence between the classical versions of natural deduction calculus (Nc), Hilbert-style deduction calculus (Hc) and Gentzen's sequent calculus (Gc).
\item Cut elimination for Gentzen's sequent calculus.
\end{enumerate}
To simplify, I only proved this for propositional calculus. Formalizing full predicate calculus gives rise to additional difficulties such as variable encodings, variable renaming, alpha-equivalence and induction over type families such as vectors.

The following table is a summary of the formalization. The number of lines in each file is rounded.\\

\begin{tabular}{c|c|c|c}
\textbf{file} & \textbf{Defines} & \textbf{Proves} & \textbf{\# lines} \\ \hline
\verb"a_base" & Variable & -- & 20  \\
\verb"b_soundness" & Provability in Nc, validity & Soundness of Nc. & 150 \\
\verb"c_completeness" & Conjunctive Normal Form & Completeness of Nc. & 400 \\
\verb"d_hilbert_calculus" & Provability in Hc & Equivalence of Nc and Hc. & 100 \\
\verb"e_sequent_calculus" & Provability in Gc & Equivalence of Nc and Gc. & 200 \\
\verb"f_cut_elimination" & Cut-free sequent calculus & Cut elimination. & 180 \\
\end{tabular}\\

In total there are just over 1000 lines, 66 theorems and 43 definitions.

In the next sections I will discuss the formalization in more detail. I will focus on the considerations which went into the formalization. I will skip all proofs and most Lemmas, but I will give informal proof sketches for the main theorems. To view all Lemmas, check the coqdoc html files. To view all proofs, check the Coq source files.

\section{Basic Definitions}\label{sec_definitions}
In \verb"a_base" I do some preparatory work, by defining variables. The relevant lines are
\begin{code}
\coqdockw{Parameter} \coqdef{a base.base mod.PropVars}{PropVars}{\coqdocaxiom{PropVars}} : \coqdockw{Set}.\coqdoceol
\coqdocnoindent
\coqdockw{Hypothesis} \coqdef{a base.base mod.Varseq dec}{Varseq\_dec}{\coqdocvariable{Varseq\_dec}} : \coqdockw{\ensuremath{\forall}} \coqdocvar{x} \coqdocvar{y}:\coqref{a base.base mod.PropVars}{\coqdocaxiom{PropVars}}, \coqexternalref{:type scope:'x7B' x 'x7D' '+' 'x7B' x 'x7D'}{http://coq.inria.fr/distrib/8.4pl3/stdlib/Coq.Init.Specif}{\coqdocnotation{\{}}\coqdocvariable{x} \coqexternalref{:type scope:x '=' x}{http://coq.inria.fr/distrib/8.4pl3/stdlib/Coq.Init.Logic}{\coqdocnotation{=}} \coqdocvariable{y}\coqexternalref{:type scope:'x7B' x 'x7D' '+' 'x7B' x 'x7D'}{http://coq.inria.fr/distrib/8.4pl3/stdlib/Coq.Init.Specif}{\coqdocnotation{\}}} \coqexternalref{:type scope:'x7B' x 'x7D' '+' 'x7B' x 'x7D'}{http://coq.inria.fr/distrib/8.4pl3/stdlib/Coq.Init.Specif}{\coqdocnotation{+}} \coqexternalref{:type scope:'x7B' x 'x7D' '+' 'x7B' x 'x7D'}{http://coq.inria.fr/distrib/8.4pl3/stdlib/Coq.Init.Specif}{\coqdocnotation{\{}}\coqdocvariable{x} \coqexternalref{:type scope:x '<>' x}{http://coq.inria.fr/distrib/8.4pl3/stdlib/Coq.Init.Logic}{\coqdocnotation{\ensuremath{\not=}}} \coqdocvariable{y}\coqexternalref{:type scope:'x7B' x 'x7D' '+' 'x7B' x 'x7D'}{http://coq.inria.fr/distrib/8.4pl3/stdlib/Coq.Init.Specif}{\coqdocnotation{\}}}.\coqdoceol
\end{code}
The first line indicates that the set of (propositional) variables is some set. The second line indicates that equality is decidable on this set. Note that since Coq is based on constructive logic, so this is not vacuous. This hypothesis is necessary for the completeness theorem. If we did not assume this, then one could take for example the constructive reals as variables. It is consistent with the constructive reals that all valuations (functions from variables to bool) are constant.\citep{bishop1967foundations} This means that it is possible that $\#0\vee\#1$ (we denote variables with `$\#$') is valid, but it is not provable. So in this case completeness would fail. These are the only axioms or parameters in the formalization. All results work even if there are only finitely many variables, even if there is just one variable, or even none at all (though that wouldn't be a very interesting case).

We can now define formulas in \verb"b_soundness".
\begin{code}
\coqdocnoindent
\coqdockw{Inductive} \coqdef{b soundness.sound mod.PropF}{PropF}{\coqdocinductive{PropF}} : \coqdockw{Set} :=\coqdoceol
\coqdocindent{0.50em}
\ensuremath{|} \coqdef{b soundness.sound mod.Var}{Var}{\coqdocconstructor{Var}} : \coqref{b soundness.X.PropVars}{\coqdocaxiom{PropVars}} \ensuremath{\rightarrow} \coqref{b soundness.PropF}{\coqdocinductive{PropF}}\coqdoceol
\coqdocindent{0.50em}
\ensuremath{|} \coqdef{b soundness.sound mod.Bot}{Bot}{\coqdocconstructor{Bot}} : \coqref{b soundness.PropF}{\coqdocinductive{PropF}}\coqdoceol
\coqdocindent{0.50em}
\ensuremath{|} \coqdef{b soundness.sound mod.Conj}{Conj}{\coqdocconstructor{Conj}} : \coqref{b soundness.PropF}{\coqdocinductive{PropF}} \ensuremath{\rightarrow} \coqref{b soundness.PropF}{\coqdocinductive{PropF}} \ensuremath{\rightarrow} \coqref{b soundness.PropF}{\coqdocinductive{PropF}}\coqdoceol
\coqdocindent{0.50em}
\ensuremath{|} \coqdef{b soundness.sound mod.Disj}{Disj}{\coqdocconstructor{Disj}} : \coqref{b soundness.PropF}{\coqdocinductive{PropF}} \ensuremath{\rightarrow} \coqref{b soundness.PropF}{\coqdocinductive{PropF}} \ensuremath{\rightarrow} \coqref{b soundness.PropF}{\coqdocinductive{PropF}}\coqdoceol
\coqdocindent{0.50em}
\ensuremath{|} \coqdef{b soundness.sound mod.Impl}{Impl}{\coqdocconstructor{Impl}} : \coqref{b soundness.PropF}{\coqdocinductive{PropF}} \ensuremath{\rightarrow} \coqref{b soundness.PropF}{\coqdocinductive{PropF}} \ensuremath{\rightarrow} \coqref{b soundness.PropF}{\coqdocinductive{PropF}}\coqdoceol
\coqdocnoindent
.\coqdoceol
\coqdocemptyline
\coqdocnoindent
\coqdockw{Notation} \coqdef{b soundness.sound mod.:My scope:'x23' x}{"}{"}\# P" := (\coqref{b soundness.sound mod.Var}{\coqdocconstructor{Var}} \coqdocvar{P}) (\coqdoctac{at} \coqdockw{level} 1) : \coqdocvar{My\_scope}.\coqdoceol
\coqdocnoindent
\coqdockw{Notation} \coqdef{b soundness.sound mod.:My scope:x 'xE2x88xA8' x}{"}{"}A ∨ B" := (\coqref{b soundness.sound mod.Disj}{\coqdocconstructor{Disj}} \coqdocvar{A} \coqdocvar{B}) (\coqdoctac{at} \coqdockw{level} 15, \coqdoctac{right} \coqdockw{associativity}) : \coqdocvar{My\_scope}.\coqdoceol
\coqdocnoindent
\coqdockw{Notation} \coqdef{b soundness.sound mod.:My scope:x 'xE2x88xA7' x}{"}{"}A ∧ B" := (\coqref{b soundness.sound mod.Conj}{\coqdocconstructor{Conj}} \coqdocvar{A} \coqdocvar{B}) (\coqdoctac{at} \coqdockw{level} 15, \coqdoctac{right} \coqdockw{associativity}) : \coqdocvar{My\_scope}.\coqdoceol
\coqdocnoindent
\coqdockw{Notation} \coqdef{b soundness.sound mod.:My scope:x 'xE2x86x92' x}{"}{"}A → B" := (\coqref{b soundness.sound mod.Impl}{\coqdocconstructor{Impl}} \coqdocvar{A} \coqdocvar{B}) (\coqdoctac{at} \coqdockw{level} 16, \coqdoctac{right} \coqdockw{associativity}) : \coqdocvar{My\_scope}.\coqdoceol
\coqdocnoindent
\coqdockw{Notation} \coqdef{b soundness.sound mod.:My scope:'xE2x8AxA5'}{"}{"}⊥" := \coqref{b soundness.sound mod.Bot}{\coqdocconstructor{Bot}} (\coqdoctac{at} \coqdockw{level} 0)  : \coqdocvar{My\_scope}.\coqdoceol
\coqdocnoindent
\coqdockw{Definition} \coqdef{b soundness.sound mod.Neg}{Neg}{\coqdocdefinition{Neg}} \coqdocvar{A} := \coqdocvariable{A} \coqref{b soundness.sound mod.:My scope:x 'xE2x86x92' x}{\coqdocnotation{→}} \coqref{b soundness.sound mod.:My scope:'xE2x8AxA5'}{\coqdocnotation{⊥}}.\coqdoceol
\coqdocnoindent
\coqdockw{Notation} \coqdef{b soundness.sound mod.:My scope:'xC2xAC' x}{"}{"}¬ A" := (\coqref{b soundness.sound mod.Neg}{\coqdocdefinition{Neg}} \coqdocvar{A}) (\coqdoctac{at} \coqdockw{level} 5) : \coqdocvar{My\_scope}.\coqdoceol
\coqdocnoindent
\coqdockw{Definition} \coqdef{b soundness.sound mod.Top}{Top}{\coqdocdefinition{Top}} := \coqref{b soundness.sound mod.:My scope:'xC2xAC' x}{\coqdocnotation{¬}}\coqref{b soundness.sound mod.:My scope:'xE2x8AxA5'}{\coqdocnotation{⊥}}.\coqdoceol
\coqdocnoindent
\coqdockw{Notation} \coqdef{b soundness.sound mod.:My scope:'xE2x8AxA4'}{"}{"}⊤" := \coqref{b soundness.sound mod.Top}{\coqdocdefinition{Top}} (\coqdoctac{at} \coqdockw{level} 0) : \coqdocvar{My\_scope}.\coqdoceol
\end{code}
In the formalization I used the unicode character $\rightarrow$ for implication, but to distinguish it from the function type of Coq, I use $\rightharpoonup$ in this document instead.

I use defined negation, because that simplifies induction proofs by having one fewer connective. Other than this, I take all connectives as primitive. I could have defined conjunction and disjunction in terms of implication and falsum and retain an equivalent system, since I only formalize classical logic. However, I decided against this, to make the formalization easier to adapt to intuitionistic or minimal logic, where these connectives aren't interdefinable. We defined negation because it is definable in all these calculi.

In an earlier version of the formalization I used a special variable $\bot$ which corresponds to the false formula. This had as advantages that the variables corresponded exactly to the atomic formulae and that there is one fewer induction base case when proving things by induction over formulae. In the end I decided against it, since it is more natural to have $\bot$ not as a variable, but as a separate constant, and the alternative required to define valuations as assigning `false' to the variable $\bot$, while no such requirement is necessary now.

A \emph{valuation} is an element of the space $\coqdocaxiom{PropVars}\to\coqdocinductive{bool}$. I define the truth of a formula $A$ under valuation $v$ in the obvious way. This allows us to define the validity of formulae. The map \coqdocdefinition{Is\_true} sends a boolean value to the corresponding (true or false) proposition. \coqdocnotation{[]} denotes the empty list.
\begin{code}
\coqdocnoindent
\coqdockw{Fixpoint} \coqdef{b soundness.sound mod.TrueQ}{TrueQ}{\coqdocdefinition{TrueQ}} \coqdocvar{v} \coqdocvar{A} : \coqexternalref{bool}{http://coq.inria.fr/distrib/8.4pl3/stdlib/Coq.Init.Datatypes}{\coqdocinductive{bool}} := \coqdockw{match} \coqdocvariable{A} \coqdockw{with}\coqdoceol
\coqdocindent{0.50em}
\ensuremath{|} \coqref{b soundness.sound mod.:My scope:'x23' x}{\coqdocnotation{\#}} \coqdocvar{P}   \ensuremath{\Rightarrow} \coqdocvariable{v} \coqdocvar{P}\coqdoceol
\coqdocindent{0.50em}
\ensuremath{|} \coqref{b soundness.sound mod.:My scope:'xE2x8AxA5'}{\coqdocnotation{⊥}}     \ensuremath{\Rightarrow} \coqexternalref{false}{http://coq.inria.fr/distrib/8.4pl3/stdlib/Coq.Init.Datatypes}{\coqdocconstructor{false}}\coqdoceol
\coqdocindent{0.50em}
\ensuremath{|} \coqdocvar{B} \coqref{b soundness.sound mod.:My scope:x 'xE2x88xA8' x}{\coqdocnotation{∨}} \coqdocvar{C} \ensuremath{\Rightarrow} \coqexternalref{:bool scope:x '||' x}{http://coq.inria.fr/distrib/8.4pl3/stdlib/Coq.Init.Datatypes}{\coqdocnotation{(}}\coqref{b soundness.TrueQ}{\coqdocdefinition{TrueQ}} \coqdocvariable{v} \coqdocvar{B}\coqexternalref{:bool scope:x '||' x}{http://coq.inria.fr/distrib/8.4pl3/stdlib/Coq.Init.Datatypes}{\coqdocnotation{)}} \coqexternalref{:bool scope:x '||' x}{http://coq.inria.fr/distrib/8.4pl3/stdlib/Coq.Init.Datatypes}{\coqdocnotation{||}} \coqexternalref{:bool scope:x '||' x}{http://coq.inria.fr/distrib/8.4pl3/stdlib/Coq.Init.Datatypes}{\coqdocnotation{(}}\coqref{b soundness.TrueQ}{\coqdocdefinition{TrueQ}} \coqdocvariable{v} \coqdocvar{C}\coqexternalref{:bool scope:x '||' x}{http://coq.inria.fr/distrib/8.4pl3/stdlib/Coq.Init.Datatypes}{\coqdocnotation{)}}\coqdoceol
\coqdocindent{0.50em}
\ensuremath{|} \coqdocvar{B} \coqref{b soundness.sound mod.:My scope:x 'xE2x88xA7' x}{\coqdocnotation{∧}} \coqdocvar{C} \ensuremath{\Rightarrow} \coqexternalref{:bool scope:x 'x26x26' x}{http://coq.inria.fr/distrib/8.4pl3/stdlib/Coq.Init.Datatypes}{\coqdocnotation{(}}\coqref{b soundness.TrueQ}{\coqdocdefinition{TrueQ}} \coqdocvariable{v} \coqdocvar{B}\coqexternalref{:bool scope:x 'x26x26' x}{http://coq.inria.fr/distrib/8.4pl3/stdlib/Coq.Init.Datatypes}{\coqdocnotation{)}} \coqexternalref{:bool scope:x 'x26x26' x}{http://coq.inria.fr/distrib/8.4pl3/stdlib/Coq.Init.Datatypes}{\coqdocnotation{\&\&}} \coqexternalref{:bool scope:x 'x26x26' x}{http://coq.inria.fr/distrib/8.4pl3/stdlib/Coq.Init.Datatypes}{\coqdocnotation{(}}\coqref{b soundness.TrueQ}{\coqdocdefinition{TrueQ}} \coqdocvariable{v} \coqdocvar{C}\coqexternalref{:bool scope:x 'x26x26' x}{http://coq.inria.fr/distrib/8.4pl3/stdlib/Coq.Init.Datatypes}{\coqdocnotation{)}}\coqdoceol
\coqdocindent{0.50em}
\ensuremath{|} \coqdocvar{B} \coqref{b soundness.sound mod.:My scope:x 'xE2x86x92' x}{\coqdocnotation{→}} \coqdocvar{C} \ensuremath{\Rightarrow} \coqexternalref{:bool scope:x '||' x}{http://coq.inria.fr/distrib/8.4pl3/stdlib/Coq.Init.Datatypes}{\coqdocnotation{(}}\coqexternalref{negb}{http://coq.inria.fr/distrib/8.4pl3/stdlib/Coq.Init.Datatypes}{\coqdocdefinition{negb}} (\coqref{b soundness.TrueQ}{\coqdocdefinition{TrueQ}} \coqdocvariable{v} \coqdocvar{B})\coqexternalref{:bool scope:x '||' x}{http://coq.inria.fr/distrib/8.4pl3/stdlib/Coq.Init.Datatypes}{\coqdocnotation{)}} \coqexternalref{:bool scope:x '||' x}{http://coq.inria.fr/distrib/8.4pl3/stdlib/Coq.Init.Datatypes}{\coqdocnotation{||}} \coqexternalref{:bool scope:x '||' x}{http://coq.inria.fr/distrib/8.4pl3/stdlib/Coq.Init.Datatypes}{\coqdocnotation{(}}\coqref{b soundness.TrueQ}{\coqdocdefinition{TrueQ}} \coqdocvariable{v} \coqdocvar{C}\coqexternalref{:bool scope:x '||' x}{http://coq.inria.fr/distrib/8.4pl3/stdlib/Coq.Init.Datatypes}{\coqdocnotation{)}}\coqdoceol
\coqdocnoindent
\coqdockw{end}.\coqdoceol
\coqdocnoindent
\coqdockw{Definition} \coqdef{b soundness.sound mod.Satisfies}{Satisfies}{\coqdocdefinition{Satisfies}} \coqdocvar{v} \coqdocvar{Γ} := \coqdockw{\ensuremath{\forall}} \coqdocvar{A}, \coqexternalref{In}{http://coq.inria.fr/distrib/8.4pl3/stdlib/Coq.Lists.List}{\coqdocdefinition{In}} \coqdocvariable{A} \coqdocvariable{Γ} \ensuremath{\rightarrow} \coqexternalref{Is true}{http://coq.inria.fr/distrib/8.4pl3/stdlib/Coq.Bool.Bool}{\coqdocdefinition{Is\_true}} (\coqref{b soundness.sound mod.TrueQ}{\coqdocdefinition{TrueQ}} \coqdocvariable{v} \coqdocvariable{A}).\coqdoceol
\coqdocnoindent
\coqdockw{Definition} \coqdef{b soundness.sound mod.Models}{Models}{\coqdocdefinition{Models}} \coqdocvar{Γ} \coqdocvar{A} := \coqdockw{\ensuremath{\forall}} \coqdocvar{v},\coqref{b soundness.sound mod.Satisfies}{\coqdocdefinition{Satisfies}} \coqdocvariable{v} \coqdocvariable{Γ}\ensuremath{\rightarrow}\coqexternalref{Is true}{http://coq.inria.fr/distrib/8.4pl3/stdlib/Coq.Bool.Bool}{\coqdocdefinition{Is\_true}} (\coqref{b soundness.sound mod.TrueQ}{\coqdocdefinition{TrueQ}} \coqdocvariable{v} \coqdocvariable{A}).\coqdoceol
\coqdocnoindent
\coqdockw{Notation} \coqdef{b soundness.sound mod.::x 'xE2x8AxA8' x}{"}{"}Γ ⊨ A" := (\coqref{b soundness.sound mod.Models}{\coqdocdefinition{Models}} \coqdocvar{Γ} \coqdocvar{A}) (\coqdoctac{at} \coqdockw{level} 80).\coqdoceol
\coqdocnoindent
\coqdockw{Definition} \coqdef{b soundness.sound mod.Valid}{Valid}{\coqdocdefinition{Valid}} \coqdocvar{A} := \coqexternalref{ListNotations.:list scope:'[' ']'}{http://coq.inria.fr/distrib/8.4pl3/stdlib/Coq.Lists.List}{\coqdocnotation{[]}} \coqref{b soundness.sound mod.::x 'xE2x8AxA8' x}{\coqdocnotation{⊨}} \coqdocvariable{A}.\coqdoceol
\end{code}
The notion of provability is naturally defined inductively. I use a context sharing version of classical natural deduction. Using the context sharing version simplifies proofs, because changes in the context are usually hard to deal with for proofs of meta-theoretic theorems. Note that I am using lists, not sets, for the context, because it is easier to reason about lists in proofs.
\begin{code}
\coqdockw{Reserved Notation} "Γ ⊢ A" (\coqdoctac{at} \coqdockw{level} 80).\coqdoceol
\coqdocnoindent
\coqdockw{Inductive} \coqdef{b soundness.sound mod.Nc}{Nc}{\coqdocinductive{Nc}} : \coqexternalref{list}{http://coq.inria.fr/distrib/8.4pl3/stdlib/Coq.Init.Datatypes}{\coqdocinductive{list}} \coqref{b soundness.sound mod.PropF}{\coqdocinductive{PropF}}\ensuremath{\rightarrow} \coqref{b soundness.sound mod.PropF}{\coqdocinductive{PropF}}\ensuremath{\rightarrow}\coqdockw{Prop} :=\coqdoceol
\coqdocnoindent
\ensuremath{|} \coqdef{b soundness.sound mod.Nax}{Nax}{\coqdocconstructor{Nax}}   : \coqdockw{\ensuremath{\forall}} \coqdocvar{Γ} \coqdocvar{A}  ,    \coqexternalref{In}{http://coq.inria.fr/distrib/8.4pl3/stdlib/Coq.Lists.List}{\coqdocdefinition{In}} \coqdocvariable{A} \coqdocvariable{Γ}                           \ensuremath{\rightarrow} \coqdocvariable{Γ} \coqref{b soundness.sound mod.:My scope:x 'xE2x8AxA2' x}{\coqdocnotation{⊢}} \coqdocvariable{A}\coqdoceol
\coqdocnoindent
\ensuremath{|} \coqdef{b soundness.sound mod.ImpI}{ImpI}{\coqdocconstructor{ImpI}}  : \coqdockw{\ensuremath{\forall}} \coqdocvar{Γ} \coqdocvar{A} \coqdocvar{B},  \coqdocvariable{A}\coqexternalref{:list scope:x '::' x}{http://coq.inria.fr/distrib/8.4pl3/stdlib/Coq.Init.Datatypes}{\coqdocnotation{::}}\coqdocvariable{Γ} \coqref{b soundness.sound mod.:My scope:x 'xE2x8AxA2' x}{\coqdocnotation{⊢}} \coqdocvariable{B}                           \ensuremath{\rightarrow} \coqdocvariable{Γ} \coqref{b soundness.sound mod.:My scope:x 'xE2x8AxA2' x}{\coqdocnotation{⊢}} \coqdocvariable{A} \coqref{b soundness.sound mod.:My scope:x 'xE2x86x92' x}{\coqdocnotation{→}} \coqdocvariable{B}\coqdoceol
\coqdocnoindent
\ensuremath{|} \coqdef{b soundness.sound mod.ImpE}{ImpE}{\coqdocconstructor{ImpE}}  : \coqdockw{\ensuremath{\forall}} \coqdocvar{Γ} \coqdocvar{A} \coqdocvar{B},     \coqdocvariable{Γ} \coqref{b soundness.sound mod.:My scope:x 'xE2x8AxA2' x}{\coqdocnotation{⊢}} \coqdocvariable{A} \coqref{b soundness.sound mod.:My scope:x 'xE2x86x92' x}{\coqdocnotation{→}} \coqdocvariable{B} \ensuremath{\rightarrow} \coqdocvariable{Γ} \coqref{b soundness.sound mod.:My scope:x 'xE2x8AxA2' x}{\coqdocnotation{⊢}} \coqdocvariable{A}              \ensuremath{\rightarrow} \coqdocvariable{Γ} \coqref{b soundness.sound mod.:My scope:x 'xE2x8AxA2' x}{\coqdocnotation{⊢}} \coqdocvariable{B}\coqdoceol
\coqdocnoindent
\ensuremath{|} \coqdef{b soundness.sound mod.BotC}{BotC}{\coqdocconstructor{BotC}}  : \coqdockw{\ensuremath{\forall}} \coqdocvar{Γ} \coqdocvar{A}  , \coqref{b soundness.sound mod.:My scope:'xC2xAC' x}{\coqdocnotation{¬}}\coqdocvariable{A}\coqexternalref{:list scope:x '::' x}{http://coq.inria.fr/distrib/8.4pl3/stdlib/Coq.Init.Datatypes}{\coqdocnotation{::}}\coqdocvariable{Γ} \coqref{b soundness.sound mod.:My scope:x 'xE2x8AxA2' x}{\coqdocnotation{⊢}} \coqref{b soundness.sound mod.:My scope:'xE2x8AxA5'}{\coqdocnotation{⊥}}                              \ensuremath{\rightarrow} \coqdocvariable{Γ} \coqref{b soundness.sound mod.:My scope:x 'xE2x8AxA2' x}{\coqdocnotation{⊢}} \coqdocvariable{A}\coqdoceol
\coqdocnoindent
\ensuremath{|} \coqdef{b soundness.sound mod.AndI}{AndI}{\coqdocconstructor{AndI}}  : \coqdockw{\ensuremath{\forall}} \coqdocvar{Γ} \coqdocvar{A} \coqdocvar{B},     \coqdocvariable{Γ} \coqref{b soundness.sound mod.:My scope:x 'xE2x8AxA2' x}{\coqdocnotation{⊢}} \coqdocvariable{A}     \ensuremath{\rightarrow} \coqdocvariable{Γ} \coqref{b soundness.sound mod.:My scope:x 'xE2x8AxA2' x}{\coqdocnotation{⊢}} \coqdocvariable{B}              \ensuremath{\rightarrow} \coqdocvariable{Γ} \coqref{b soundness.sound mod.:My scope:x 'xE2x8AxA2' x}{\coqdocnotation{⊢}} \coqdocvariable{A}\coqref{b soundness.sound mod.:My scope:x 'xE2x88xA7' x}{\coqdocnotation{∧}}\coqdocvariable{B}\coqdoceol
\coqdocnoindent
\ensuremath{|} \coqdef{b soundness.sound mod.AndE1}{AndE1}{\coqdocconstructor{AndE1}} : \coqdockw{\ensuremath{\forall}} \coqdocvar{Γ} \coqdocvar{A} \coqdocvar{B},     \coqdocvariable{Γ} \coqref{b soundness.sound mod.:My scope:x 'xE2x8AxA2' x}{\coqdocnotation{⊢}} \coqdocvariable{A}\coqref{b soundness.sound mod.:My scope:x 'xE2x88xA7' x}{\coqdocnotation{∧}}\coqdocvariable{B}                        \ensuremath{\rightarrow} \coqdocvariable{Γ} \coqref{b soundness.sound mod.:My scope:x 'xE2x8AxA2' x}{\coqdocnotation{⊢}} \coqdocvariable{A}\coqdoceol
\coqdocnoindent
\ensuremath{|} \coqdef{b soundness.sound mod.AndE2}{AndE2}{\coqdocconstructor{AndE2}} : \coqdockw{\ensuremath{\forall}} \coqdocvar{Γ} \coqdocvar{A} \coqdocvar{B},     \coqdocvariable{Γ} \coqref{b soundness.sound mod.:My scope:x 'xE2x8AxA2' x}{\coqdocnotation{⊢}} \coqdocvariable{A}\coqref{b soundness.sound mod.:My scope:x 'xE2x88xA7' x}{\coqdocnotation{∧}}\coqdocvariable{B}                        \ensuremath{\rightarrow} \coqdocvariable{Γ} \coqref{b soundness.sound mod.:My scope:x 'xE2x8AxA2' x}{\coqdocnotation{⊢}} \coqdocvariable{B}\coqdoceol
\coqdocnoindent
\ensuremath{|} \coqdef{b soundness.sound mod.OrI1}{OrI1}{\coqdocconstructor{OrI1}}  : \coqdockw{\ensuremath{\forall}} \coqdocvar{Γ} \coqdocvar{A} \coqdocvar{B},     \coqdocvariable{Γ} \coqref{b soundness.sound mod.:My scope:x 'xE2x8AxA2' x}{\coqdocnotation{⊢}} \coqdocvariable{A}                           \ensuremath{\rightarrow} \coqdocvariable{Γ} \coqref{b soundness.sound mod.:My scope:x 'xE2x8AxA2' x}{\coqdocnotation{⊢}} \coqdocvariable{A}\coqref{b soundness.sound mod.:My scope:x 'xE2x88xA8' x}{\coqdocnotation{∨}}\coqdocvariable{B}\coqdoceol
\coqdocnoindent
\ensuremath{|} \coqdef{b soundness.sound mod.OrI2}{OrI2}{\coqdocconstructor{OrI2}}  : \coqdockw{\ensuremath{\forall}} \coqdocvar{Γ} \coqdocvar{A} \coqdocvar{B},     \coqdocvariable{Γ} \coqref{b soundness.sound mod.:My scope:x 'xE2x8AxA2' x}{\coqdocnotation{⊢}} \coqdocvariable{B}                           \ensuremath{\rightarrow} \coqdocvariable{Γ} \coqref{b soundness.sound mod.:My scope:x 'xE2x8AxA2' x}{\coqdocnotation{⊢}} \coqdocvariable{A}\coqref{b soundness.sound mod.:My scope:x 'xE2x88xA8' x}{\coqdocnotation{∨}}\coqdocvariable{B}\coqdoceol
\coqdocnoindent
\ensuremath{|} \coqdef{b soundness.sound mod.OrE}{OrE}{\coqdocconstructor{OrE}}   : \coqdockw{\ensuremath{\forall}} \coqdocvar{Γ} \coqdocvar{A} \coqdocvar{B} \coqdocvar{C},   \coqdocvariable{Γ} \coqref{b soundness.sound mod.:My scope:x 'xE2x8AxA2' x}{\coqdocnotation{⊢}} \coqdocvariable{A}\coqref{b soundness.sound mod.:My scope:x 'xE2x88xA8' x}{\coqdocnotation{∨}}\coqdocvariable{B} \ensuremath{\rightarrow} \coqdocvariable{A}\coqexternalref{:list scope:x '::' x}{http://coq.inria.fr/distrib/8.4pl3/stdlib/Coq.Init.Datatypes}{\coqdocnotation{::}}\coqdocvariable{Γ} \coqref{b soundness.sound mod.:My scope:x 'xE2x8AxA2' x}{\coqdocnotation{⊢}} \coqdocvariable{C} \ensuremath{\rightarrow} \coqdocvariable{B}\coqexternalref{:list scope:x '::' x}{http://coq.inria.fr/distrib/8.4pl3/stdlib/Coq.Init.Datatypes}{\coqdocnotation{::}}\coqdocvariable{Γ} \coqref{b soundness.sound mod.:My scope:x 'xE2x8AxA2' x}{\coqdocnotation{⊢}} \coqdocvariable{C} \ensuremath{\rightarrow} \coqdocvariable{Γ} \coqref{b soundness.sound mod.:My scope:x 'xE2x8AxA2' x}{\coqdocnotation{⊢}} \coqdocvariable{C}\coqdoceol
\coqdocnoindent
\coqdockw{where} \coqdef{b soundness.sound mod.:My scope:x 'xE2x8AxA2' x}{"}{"}Γ ⊢ A" := (\coqref{b soundness.Nc}{\coqdocinductive{Nc}} \coqdocvar{Γ} \coqdocvar{A}) : \coqdocvar{My\_scope}.\coqdoceol
\coqdocemptyline
\coqdocnoindent
\coqdockw{Definition} \coqdef{b soundness.sound mod.Provable}{Provable}{\coqdocdefinition{Provable}} \coqdocvar{A} := \coqexternalref{ListNotations.:list scope:'[' ']'}{http://coq.inria.fr/distrib/8.4pl3/stdlib/Coq.Lists.List}{\coqdocnotation{[]}} \coqref{b soundness.sound mod.:My scope:x 'xE2x8AxA2' x}{\coqdocnotation{⊢}} \coqdocvariable{A}.\coqdoceol
\end{code}
This allows us to define the propositions we're aiming to prove.
\begin{code}
\coqdockw{Definition} \coqdef{b soundness.sound mod.Prop Soundness}{Prop\_Soundness}{\coqdocdefinition{Prop\_Soundness}} := \coqdockw{\ensuremath{\forall}} \coqdocvar{A},\coqref{b soundness.sound mod.Provable}{\coqdocdefinition{Provable}} \coqdocvariable{A}\ensuremath{\rightarrow}\coqref{b soundness.sound mod.Valid}{\coqdocdefinition{Valid}} \coqdocvariable{A}.\coqdoceol
\coqdocnoindent
\coqdockw{Definition} \coqdef{b soundness.sound mod.Prop Completeness}{Prop\_Completeness}{\coqdocdefinition{Prop\_Completeness}} := \coqdockw{\ensuremath{\forall}} \coqdocvar{A},\coqref{b soundness.sound mod.Valid}{\coqdocdefinition{Valid}} \coqdocvariable{A}\ensuremath{\rightarrow}\coqref{b soundness.sound mod.Provable}{\coqdocdefinition{Provable}} \coqdocvariable{A}.\coqdoceol
\end{code}
Given those definitions, Soundness can be proved directly by induction on the derivation (no additional lemmas are needed). Indeed, in the formalization this only requires 17 lines.\footnote{Counting the lines for theorems \coqdoclemma{Soundness\_general} and \coqdoclemma{Soundness} and for tactics \coqdocvar{case\_bool} and \coqdocvar{prove\_satisfaction}. No other lemmas or tactics are used.}
\begin{code}
\coqdocnoindent
\coqdockw{Theorem} \coqdef{b soundness.sound mod.Soundness}{Soundness}{\coqdoclemma{Soundness}} : \coqref{b soundness.sound mod.Prop Soundness}{\coqdocdefinition{Prop\_Soundness}}.\coqdoceol
\end{code}
\section{Completeness}\label{sec_complete}

The formalization continues to prove Completeness. I will first give a proof sketch, before I turn my attention to the formalization.

\subsection{Proof Sketch}\label{sec_complete_sketch}

The formalization uses the following proof. Recall that a literal is an atomic formula or the negation of an atomic formula, a clause is a disjunction of literals, and that a formula is in conjunctive normal form (CNF) if it's a conjunction of clauses. Also, a formula is in \emph{negation normal form} (NNF), if it only consists of conjunctions, disjunctions and literals. In this paper I will use CNF and NNF as nouns indicating formulae in CNF resp. NNF.\\

\emph{Proof sketch of Completeness}. Let $A$ be a formula, $A\nnf$ its negation normal form and $A\cnf$ its conjunctive normal form. We say that a clause is \emph{syntactically invalid} if either it contains both $p$ and $\neg p$ for some atomic formula $p$, or if it contains $\neg\bot$. We say that $A\cnf$ is \emph{syntactically valid} if it contains no syntactically invalid clauses.
Then Completeness follows from the following statements
\begin{enumerate}[(1)]
\item If $A$ is valid, then $A\nnf$ is valid.
\item If $A\nnf$ is valid, then $A\cnf$ is valid.
\item If $A\cnf$ is valid then $A\cnf$ is syntactically valid.
\item If $A\cnf$ is syntactically valid then $A\cnf$ is provable.
\item $A\cnf\to A\nnf$ is provable.
\item $A\nnf\to A$ is provable.
\end{enumerate}
All statements except (3) are proven by induction to the structure of the formula in question. Many of these proofs need additional lemmas, which are also proven by induction. Statement (3) is most naturally proved using contraposition, by showing that if $A\cnf$ is syntactically invalid, then there is a countervaluation $v$ for it. First one does this for a clause. The countervaluation $v$ of a syntactically invalid clause $C$ can be defined as $v(x)=\text{true}$ iff $\neg x$ occurs in $C$. One can then show that $v$ is indeed a valuation, and that $v$ makes $C$ false. Then if $A\cnf$ is syntactically invalid, it contains some syntactically invalid clause $C$, and the countervaluation for that clause is also a countervaluation for $A\cnf$.\qed\\

\subsection{Considerations for the formalization}\label{sec_complete_considerations}
In the formalization I defined clauses as lists of literals and CNFs (formulae in CNF) as lists of clauses. The alternative, to say that only some class of propositional formulae are in CNF, will be very hard to work with. One of the problems which one has to deal with is associativity. Which of the following do you call clauses?
$$(((x\vee\neg y)\vee z)\vee w\qquad (x\vee(\neg y\vee (z\vee w)))\qquad (x\vee\neg y)\vee (z\vee w).$$
You can call all of them clauses, but then clauses are not in a canonical form, and proving theorems about them will probably be harder. Also, induction proofs will be harder, because you're constantly proving or using the fact that some formula is in the shape of a CNF. Using lists, you can make sure that all clauses have the same shape of parentheses, and a term of the correct type is automatically in CNF. The disadvantage, that you have to define a mapping from clauses to actual formulas, is very minor.

A viable alternative to using lists would be using non-empty lists, i.e. lists which have at least one element. When transforming a formula to conjunctive normal form, it turns out we never end up with an empty clause (i.e. an empty list of literals) or an empty formula in CNF (i.e. an empty list of clauses). Using non-empty lists might be more natural in that case, but I decided against this, because I thought this would make meta-theorems harder to prove. With non-empty lists, the base case is a list of one element, and that case is probably harder than the case of an empty list. However, allowing for empty lists gives an unexpected artefact in the translation from lists to formulae. When transforming a clause (a list of literals) to a formulae, it is natural to say $f(\coqdocconstructor{nil})=\bot$ (\coqdocconstructor{nil} is the empty list) and $f(a\coqdocnotation{::}l)=a\vee f(l)$. However, this means that, for example, $f([x;y;z])=x\vee (y\vee(z\vee\bot))$ instead of the more natural $f([x;y;z])=x\vee (y\vee z)$. This poses no problems except that it is unnatural.

\subsection{Definitions}\label{sec_complete_defs}

The definitions needed for Completeness are given below, given in \verb"c_completeness".
\begin{code}
\coqdocnoindent
\coqdockw{Inductive} \coqdef{c completeness.complete mod.NNF}{NNF}{\coqdocinductive{NNF}} : \coqdockw{Set} :=\coqdoceol
\coqdocindent{0.50em}
\ensuremath{|} \coqdef{c completeness.complete mod.NPos}{NPos}{\coqdocconstructor{NPos}} : \coqref{c completeness.X.PropVars}{\coqdocaxiom{PropVars}} \ensuremath{\rightarrow} \coqref{c completeness.NNF}{\coqdocinductive{NNF}}\coqdoceol
\coqdocindent{0.50em}
\ensuremath{|} \coqdef{c completeness.complete mod.NNeg}{NNeg}{\coqdocconstructor{NNeg}} : \coqref{c completeness.X.PropVars}{\coqdocaxiom{PropVars}} \ensuremath{\rightarrow} \coqref{c completeness.NNF}{\coqdocinductive{NNF}}\coqdoceol
\coqdocindent{0.50em}
\ensuremath{|} \coqdef{c completeness.complete mod.NBot}{NBot}{\coqdocconstructor{NBot}} : \coqref{c completeness.NNF}{\coqdocinductive{NNF}}\coqdoceol
\coqdocindent{0.50em}
\ensuremath{|} \coqdef{c completeness.complete mod.NTop}{NTop}{\coqdocconstructor{NTop}} : \coqref{c completeness.NNF}{\coqdocinductive{NNF}}\coqdoceol
\coqdocindent{0.50em}
\ensuremath{|} \coqdef{c completeness.complete mod.NConj}{NConj}{\coqdocconstructor{NConj}} : \coqref{c completeness.NNF}{\coqdocinductive{NNF}} \ensuremath{\rightarrow} \coqref{c completeness.NNF}{\coqdocinductive{NNF}} \ensuremath{\rightarrow} \coqref{c completeness.NNF}{\coqdocinductive{NNF}}\coqdoceol
\coqdocindent{0.50em}
\ensuremath{|} \coqdef{c completeness.complete mod.NDisj}{NDisj}{\coqdocconstructor{NDisj}} : \coqref{c completeness.NNF}{\coqdocinductive{NNF}} \ensuremath{\rightarrow} \coqref{c completeness.NNF}{\coqdocinductive{NNF}} \ensuremath{\rightarrow} \coqref{c completeness.NNF}{\coqdocinductive{NNF}}\coqdoceol
\coqdocnoindent
\coqdockw{Fixpoint} \coqdef{c completeness.complete mod.MakeNNF}{MakeNNF}{\coqdocdefinition{MakeNNF}} (\coqdocvar{A}:\coqref{c completeness.Y.PropF}{\coqdocinductive{PropF}}) : \coqref{c completeness.complete mod.NNF}{\coqdocinductive{NNF}} := \coqdockw{match} \coqdocvariable{A} \coqdockw{with}\coqdoceol
\coqdocindent{0.50em}
\ensuremath{|} \coqref{b soundness.sound mod.:My scope:'x23' x}{\coqdocnotation{\#}} \coqdocvar{P}   \ensuremath{\Rightarrow} \coqref{c completeness.complete mod.NPos}{\coqdocconstructor{NPos}} \coqdocvar{P}\coqdoceol
\coqdocindent{0.50em}
\ensuremath{|} \coqref{b soundness.sound mod.:My scope:'xE2x8AxA5'}{\coqdocnotation{⊥}}     \ensuremath{\Rightarrow} \coqref{c completeness.complete mod.NBot}{\coqdocconstructor{NBot}}\coqdoceol
\coqdocindent{0.50em}
\ensuremath{|} \coqdocvar{B} \coqref{b soundness.sound mod.:My scope:x 'xE2x88xA8' x}{\coqdocnotation{∨}} \coqdocvar{C} \ensuremath{\Rightarrow} \coqref{c completeness.complete mod.NDisj}{\coqdocconstructor{NDisj}} (\coqref{c completeness.MakeNNF}{\coqdocdefinition{MakeNNF}} \coqdocvar{B}) (\coqref{c completeness.MakeNNF}{\coqdocdefinition{MakeNNF}} \coqdocvar{C})\coqdoceol
\coqdocindent{0.50em}
\ensuremath{|} \coqdocvar{B} \coqref{b soundness.sound mod.:My scope:x 'xE2x88xA7' x}{\coqdocnotation{∧}} \coqdocvar{C} \ensuremath{\Rightarrow} \coqref{c completeness.complete mod.NConj}{\coqdocconstructor{NConj}} (\coqref{c completeness.MakeNNF}{\coqdocdefinition{MakeNNF}} \coqdocvar{B}) (\coqref{c completeness.MakeNNF}{\coqdocdefinition{MakeNNF}} \coqdocvar{C})\coqdoceol
\coqdocindent{0.50em}
\ensuremath{|} \coqdocvar{B} \coqref{b soundness.sound mod.:My scope:x 'xE2x86x92' x}{\coqdocnotation{→}} \coqdocvar{C} \ensuremath{\Rightarrow} \coqref{c completeness.complete mod.NDisj}{\coqdocconstructor{NDisj}} (\coqref{c completeness.MakeNNFN}{\coqdocdefinition{MakeNNFN}} \coqdocvar{B}) (\coqref{c completeness.MakeNNF}{\coqdocdefinition{MakeNNF}} \coqdocvar{C})\coqdoceol
\coqdocindent{0.50em}
\coqdockw{end}\coqdoceol
\coqdocnoindent
\coqdockw{with} \coqdef{c completeness.complete mod.MakeNNFN}{MakeNNFN}{\coqdocdefinition{MakeNNFN}} (\coqdocvar{A}:\coqref{c completeness.Y.PropF}{\coqdocinductive{PropF}}) : \coqref{c completeness.complete mod.NNF}{\coqdocinductive{NNF}} := \coqdockw{match} \coqdocvariable{A} \coqdockw{with}\coqdoceol
\coqdocindent{0.50em}
\ensuremath{|} \coqref{b soundness.sound mod.:My scope:'x23' x}{\coqdocnotation{\#}} \coqdocvar{P}   \ensuremath{\Rightarrow} \coqref{c completeness.complete mod.NNeg}{\coqdocconstructor{NNeg}} \coqdocvar{P}\coqdoceol
\coqdocindent{0.50em}
\ensuremath{|} \coqref{b soundness.sound mod.:My scope:'xE2x8AxA5'}{\coqdocnotation{⊥}}     \ensuremath{\Rightarrow} \coqref{c completeness.complete mod.NTop}{\coqdocconstructor{NTop}}\coqdoceol
\coqdocindent{0.50em}
\ensuremath{|} \coqdocvar{B} \coqref{b soundness.sound mod.:My scope:x 'xE2x88xA8' x}{\coqdocnotation{∨}} \coqdocvar{C} \ensuremath{\Rightarrow} \coqref{c completeness.complete mod.NConj}{\coqdocconstructor{NConj}} (\coqref{c completeness.MakeNNFN}{\coqdocdefinition{MakeNNFN}} \coqdocvar{B}) (\coqref{c completeness.MakeNNFN}{\coqdocdefinition{MakeNNFN}} \coqdocvar{C})\coqdoceol
\coqdocindent{0.50em}
\ensuremath{|} \coqdocvar{B} \coqref{b soundness.sound mod.:My scope:x 'xE2x88xA7' x}{\coqdocnotation{∧}} \coqdocvar{C} \ensuremath{\Rightarrow} \coqref{c completeness.complete mod.NDisj}{\coqdocconstructor{NDisj}} (\coqref{c completeness.MakeNNFN}{\coqdocdefinition{MakeNNFN}} \coqdocvar{B}) (\coqref{c completeness.MakeNNFN}{\coqdocdefinition{MakeNNFN}} \coqdocvar{C})\coqdoceol
\coqdocindent{0.50em}
\ensuremath{|} \coqdocvar{B} \coqref{b soundness.sound mod.:My scope:x 'xE2x86x92' x}{\coqdocnotation{→}} \coqdocvar{C} \ensuremath{\Rightarrow} \coqref{c completeness.complete mod.NConj}{\coqdocconstructor{NConj}} (\coqref{c completeness.MakeNNF}{\coqdocdefinition{MakeNNF}} \coqdocvar{B}) (\coqref{c completeness.MakeNNFN}{\coqdocdefinition{MakeNNFN}} \coqdocvar{C})\coqdoceol
\coqdocindent{0.50em}
\coqdockw{end}.\coqdoceol
\coqdocnoindent
\coqdockw{Inductive} \coqdef{c completeness.complete mod.Literal}{Literal}{\coqdocinductive{Literal}} :=\coqdoceol
\coqdocnoindent
\ensuremath{|} \coqdef{c completeness.complete mod.LPos}{LPos}{\coqdocconstructor{LPos}} : \coqref{c completeness.X.PropVars}{\coqdocaxiom{PropVars}} \ensuremath{\rightarrow} \coqref{c completeness.Literal}{\coqdocinductive{Literal}}\coqdoceol
\coqdocnoindent
\ensuremath{|} \coqdef{c completeness.complete mod.LNeg}{LNeg}{\coqdocconstructor{LNeg}} : \coqref{c completeness.X.PropVars}{\coqdocaxiom{PropVars}} \ensuremath{\rightarrow} \coqref{c completeness.Literal}{\coqdocinductive{Literal}}\coqdoceol
\coqdocnoindent
\ensuremath{|} \coqdef{c completeness.complete mod.LBot}{LBot}{\coqdocconstructor{LBot}} : \coqref{c completeness.Literal}{\coqdocinductive{Literal}}\coqdoceol
\coqdocnoindent
\ensuremath{|} \coqdef{c completeness.complete mod.LTop}{LTop}{\coqdocconstructor{LTop}} : \coqref{c completeness.Literal}{\coqdocinductive{Literal}}\coqdoceol
\coqdocnoindent
.\coqdoceol
\end{code}
The inclusions \coqdocdefinition{NNFtoPropF} and \coqdocdefinition{LiteraltoPropF} from respectively NNF and literals to propositional formulae are defined in the obvious way.

In \verb"a_base" I define the constructor \coqdocdefinition{map\_fold\_right}, such that for $f : B\to A$, $g : A\to  A\to  A$, $a:A$ and $x_i:B$ we have
$$\coqdocdefinition{map\_fold\_right}\,f\,g\,a\,[x_0; \cdots ; x_n] = g (f(x_0),g (f(x_1),\cdots g (f(x_n),a)\cdots).$$
The definition is below, and allows us to define the inclusion from Clause and CNF to PropF.
\begin{code}
\coqdocnoindent
\coqdockw{Fixpoint} \coqdef{a base.base mod.map fold right}{map\_fold\_right}{\coqdocdefinition{map\_fold\_right}} (\coqdocvar{A} \coqdocvar{B}:\coqdockw{Type}) (\coqdocvar{f} : \coqdocvariable{B} \ensuremath{\rightarrow} \coqdocvariable{A}) (\coqdocvar{g} : \coqdocvariable{A} \ensuremath{\rightarrow} \coqdocvariable{A} \ensuremath{\rightarrow} \coqdocvariable{A}) \coqdocvar{a} \coqdocvar{l} := \coqdockw{match} \coqdocvariable{l} \coqdockw{with}\coqdoceol
\coqdocindent{0.50em}
\ensuremath{|} \coqexternalref{nil}{http://coq.inria.fr/distrib/8.4pl3/stdlib/Coq.Init.Datatypes}{\coqdocconstructor{nil}} \ensuremath{\Rightarrow} \coqdocvariable{a}\coqdoceol
\coqdocindent{0.50em}
\ensuremath{|} \coqdocvar{b}\coqexternalref{:list scope:x '::' x}{http://coq.inria.fr/distrib/8.4pl3/stdlib/Coq.Init.Datatypes}{\coqdocnotation{::}}\coqdocvar{l2} \ensuremath{\Rightarrow} \coqdocvariable{g} (\coqdocvariable{f} \coqdocvar{b}) (\coqref{a base.map fold right}{\coqdocdefinition{map\_fold\_right}} \coqdocvariable{f} \coqdocvariable{g} \coqdocvariable{a} \coqdocvar{l2})\coqdoceol
\coqdocnoindent
\coqdockw{end}.\coqdoceol
\coqdocnoindent
\coqdockw{Definition} \coqdef{c completeness.complete mod.Clause}{Clause}{\coqdocdefinition{Clause}} := \coqexternalref{list}{http://coq.inria.fr/distrib/8.4pl3/stdlib/Coq.Init.Datatypes}{\coqdocinductive{list}} \coqref{c completeness.complete mod.Literal}{\coqdocinductive{Literal}}.\coqdoceol
\coqdocnoindent
\coqdockw{Definition} \coqdef{c completeness.complete mod.ClausetoPropF}{ClausetoPropF}{\coqdocdefinition{ClausetoPropF}} := \coqref{c completeness.X.map fold right}{\coqdocdefinition{map\_fold\_right}} \coqref{c completeness.complete mod.LiteraltoPropF}{\coqdocdefinition{LiteraltoPropF}} \coqref{c completeness.Y.Disj}{\coqdocconstructor{Disj}} \coqref{b soundness.sound mod.:My scope:'xE2x8AxA5'}{\coqdocnotation{⊥}}.\coqdoceol
\coqdocemptyline
\coqdocnoindent
\coqdockw{Definition} \coqdef{c completeness.complete mod.CNF}{CNF}{\coqdocdefinition{CNF}} := \coqexternalref{list}{http://coq.inria.fr/distrib/8.4pl3/stdlib/Coq.Init.Datatypes}{\coqdocinductive{list}} \coqref{c completeness.complete mod.Clause}{\coqdocdefinition{Clause}}.\coqdoceol
\coqdocemptyline
\coqdocnoindent
\coqdockw{Definition} \coqdef{c completeness.complete mod.CNFtoPropF}{CNFtoPropF}{\coqdocdefinition{CNFtoPropF}} := \coqref{c completeness.X.map fold right}{\coqdocdefinition{map\_fold\_right}} \coqref{c completeness.complete mod.ClausetoPropF}{\coqdocdefinition{ClausetoPropF}} \coqref{c completeness.Y.Conj}{\coqdocconstructor{Conj}} \coqref{b soundness.sound mod.:My scope:'xE2x8AxA4'}{\coqdocnotation{⊤}}.\coqdoceol
\end{code}
We still need to define the map which transforms a NNF to a CNF. The conjunction of two CNFs is just concatenation of the corresponding lists. However, the disjunction is harder to define. Doing first the simpler case of taking the disjunction of a clause with a CNF, we see that this corresponds to just adding this clause in front of every other clause:
$$C\vee(C_1\wedge C_2\wedge\cdots)\rightsquigarrow (C\vee C_1)\wedge(C\vee C_2)\wedge\cdots,$$
see \coqdocdefinition{AddClause} below ($\coqdocdefinition{map}\,f\,l$ means applying $f$ to every element of $l$). We can then define the disjunction of two CNFs as follows:
$$(C_1\wedge C_2\wedge\cdots)\vee A\rightsquigarrow (C_1\vee A)\wedge (C_2\vee A)\wedge\cdots$$
where $A$ is an CNF, and where $C_1\vee A$ is defined using \coqdocdefinition{AddClause}. This gives rise to the definition of \coqdocdefinition{Disjunct} below. The term \coqdocdefinition{flat\_map} is defined in the Coq library such that if $f:A\to\coqdocinductive{list} B$ and $x_i:A$ then
$$\coqdocdefinition{flat\_map}\,f\,[x_0; \cdots;x_n] = f (x_0)\;\coqdocnotation{++}\; f(x_1)\;\coqdocnotation{++}\;\cdots\; \coqdocnotation{++}\;f(x_n).$$
This allows us to define the transformation from NNF to CNF.
\begin{code}
\coqdocnoindent
\coqdockw{Definition} \coqdef{c completeness.complete mod.AddClause}{AddClause}{\coqdocdefinition{AddClause}} (\coqdocvar{l}:\coqref{c completeness.complete mod.Clause}{\coqdocdefinition{Clause}}) (\coqdocvar{ll}:\coqref{c completeness.complete mod.CNF}{\coqdocdefinition{CNF}}) : \coqref{c completeness.complete mod.CNF}{\coqdocdefinition{CNF}} := \coqexternalref{map}{http://coq.inria.fr/distrib/8.4pl3/stdlib/Coq.Lists.List}{\coqdocdefinition{map}} (\coqdockw{fun} \coqdocvar{l2} \ensuremath{\Rightarrow} \coqdocvariable{l}\coqexternalref{:list scope:x '++' x}{http://coq.inria.fr/distrib/8.4pl3/stdlib/Coq.Init.Datatypes}{\coqdocnotation{++}}\coqdocvariable{l2}) \coqdocvariable{ll}.\coqdoceol
\coqdocnoindent
\coqdockw{Definition} \coqdef{c completeness.complete mod.Disjunct}{Disjunct}{\coqdocdefinition{Disjunct}} (\coqdocvar{ll} \coqdocvar{ll2}:\coqref{c completeness.complete mod.CNF}{\coqdocdefinition{CNF}}) : \coqref{c completeness.complete mod.CNF}{\coqdocdefinition{CNF}} := \coqexternalref{flat map}{http://coq.inria.fr/distrib/8.4pl3/stdlib/Coq.Lists.List}{\coqdocdefinition{flat\_map}} (\coqdockw{fun} \coqdocvar{l} \ensuremath{\Rightarrow} \coqref{c completeness.complete mod.AddClause}{\coqdocdefinition{AddClause}} \coqdocvariable{l} \coqdocvariable{ll2}) \coqdocvariable{ll}.\coqdoceol
\coqdocemptyline
\coqdocnoindent
\coqdockw{Fixpoint} \coqdef{c completeness.complete mod.MakeCNF}{MakeCNF}{\coqdocdefinition{MakeCNF}} (\coqdocvar{A}:\coqref{c completeness.complete mod.NNF}{\coqdocinductive{NNF}}) : \coqref{c completeness.complete mod.CNF}{\coqdocdefinition{CNF}} := \coqdockw{match} \coqdocvariable{A} \coqdockw{with}\coqdoceol
\coqdocindent{0.50em}
\ensuremath{|} \coqref{c completeness.complete mod.NPos}{\coqdocconstructor{NPos}} \coqdocvar{P}    \ensuremath{\Rightarrow} \coqexternalref{ListNotations.:list scope:'[' x ';' '..' ';' x ']'}{http://coq.inria.fr/distrib/8.4pl3/stdlib/Coq.Lists.List}{\coqdocnotation{[[}}\coqref{c completeness.complete mod.LPos}{\coqdocconstructor{LPos}} \coqdocvar{P}\coqexternalref{ListNotations.:list scope:'[' x ';' '..' ';' x ']'}{http://coq.inria.fr/distrib/8.4pl3/stdlib/Coq.Lists.List}{\coqdocnotation{]]}}\coqdoceol
\coqdocindent{0.50em} \ensuremath{|} \coqref{c completeness.complete mod.NNeg}{\coqdocconstructor{NNeg}} \coqdocvar{P}    \ensuremath{\Rightarrow} \coqexternalref{ListNotations.:list scope:'[' x ';' '..' ';' x ']'}{http://coq.inria.fr/distrib/8.4pl3/stdlib/Coq.Lists.List}{\coqdocnotation{[[}}\coqref{c completeness.complete mod.LNeg}{\coqdocconstructor{LNeg}} \coqdocvar{P}\coqexternalref{ListNotations.:list scope:'[' x ';' '..' ';' x ']'}{http://coq.inria.fr/distrib/8.4pl3/stdlib/Coq.Lists.List}{\coqdocnotation{]]}}\coqdoceol
\coqdocindent{0.50em} \ensuremath{|} \coqref{c completeness.complete mod.NBot}{\coqdocconstructor{NBot}}      \ensuremath{\Rightarrow} \coqexternalref{ListNotations.:list scope:'[' x ';' '..' ';' x ']'}{http://coq.inria.fr/distrib/8.4pl3/stdlib/Coq.Lists.List}{\coqdocnotation{[[}}\coqref{c completeness.complete mod.LBot}{\coqdocconstructor{LBot}}\coqexternalref{ListNotations.:list scope:'[' x ';' '..' ';' x ']'}{http://coq.inria.fr/distrib/8.4pl3/stdlib/Coq.Lists.List}{\coqdocnotation{]]}}\coqdoceol
\coqdocindent{0.50em} \ensuremath{|} \coqref{c completeness.complete mod.NTop}{\coqdocconstructor{NTop}}      \ensuremath{\Rightarrow} \coqexternalref{ListNotations.:list scope:'[' x ';' '..' ';' x ']'}{http://coq.inria.fr/distrib/8.4pl3/stdlib/Coq.Lists.List}{\coqdocnotation{[[}}\coqref{c completeness.complete mod.LTop}{\coqdocconstructor{LTop}}\coqexternalref{ListNotations.:list scope:'[' x ';' '..' ';' x ']'}{http://coq.inria.fr/distrib/8.4pl3/stdlib/Coq.Lists.List}{\coqdocnotation{]]}}\coqdoceol
\coqdocindent{0.50em} \ensuremath{|} \coqref{c completeness.complete mod.NConj}{\coqdocconstructor{NConj}} \coqdocvar{B} \coqdocvar{C} \ensuremath{\Rightarrow} \coqref{c completeness.MakeCNF}{\coqdocdefinition{MakeCNF}} \coqdocvar{B} \coqexternalref{:list scope:x '++' x}{http://coq.inria.fr/distrib/8.4pl3/stdlib/Coq.Init.Datatypes}{\coqdocnotation{++}} \coqref{c completeness.MakeCNF}{\coqdocdefinition{MakeCNF}} \coqdocvar{C}\coqdoceol
\coqdocindent{0.50em}
\ensuremath{|} \coqref{c completeness.complete mod.NDisj}{\coqdocconstructor{NDisj}} \coqdocvar{B} \coqdocvar{C} \ensuremath{\Rightarrow} \coqref{c completeness.complete mod.Disjunct}{\coqdocdefinition{Disjunct}} (\coqref{c completeness.MakeCNF}{\coqdocdefinition{MakeCNF}} \coqdocvar{B}) (\coqref{c completeness.MakeCNF}{\coqdocdefinition{MakeCNF}} \coqdocvar{C})\coqdoceol
\coqdocindent{0.50em}
\coqdockw{end}.\coqdoceol
\end{code}
Finally, we can define the syntactical validity.
\begin{code}
  \coqdocnoindent
\coqdockw{Definition} \coqdef{c completeness.complete mod.Valid Clause}{Valid\_Clause}{\coqdocdefinition{Valid\_Clause}} (\coqdocvar{l}:\coqref{c completeness.complete mod.Clause}{\coqdocdefinition{Clause}}) := \coqexternalref{In}{http://coq.inria.fr/distrib/8.4pl3/stdlib/Coq.Lists.List}{\coqdocdefinition{In}} \coqref{c completeness.complete mod.LTop}{\coqdocconstructor{LTop}} \coqdocvariable{l}\coqexternalref{:type scope:x 'x5C/' x}{http://coq.inria.fr/distrib/8.4pl3/stdlib/Coq.Init.Logic}{\coqdocnotation{\ensuremath{\lor}}}\coqexternalref{:type scope:'exists' x '..' x ',' x}{http://coq.inria.fr/distrib/8.4pl3/stdlib/Coq.Init.Logic}{\coqdocnotation{\ensuremath{\exists}}} \coqdocvar{A}\coqexternalref{:type scope:'exists' x '..' x ',' x}{http://coq.inria.fr/distrib/8.4pl3/stdlib/Coq.Init.Logic}{\coqdocnotation{,(}}\coqexternalref{In}{http://coq.inria.fr/distrib/8.4pl3/stdlib/Coq.Lists.List}{\coqdocdefinition{In}} (\coqref{c completeness.complete mod.LPos}{\coqdocconstructor{LPos}} \coqdocvariable{A}) \coqdocvariable{l}\coqexternalref{:type scope:x '/x5C' x}{http://coq.inria.fr/distrib/8.4pl3/stdlib/Coq.Init.Logic}{\coqdocnotation{\ensuremath{\land}}}\coqexternalref{In}{http://coq.inria.fr/distrib/8.4pl3/stdlib/Coq.Lists.List}{\coqdocdefinition{In}} (\coqref{c completeness.complete mod.LNeg}{\coqdocconstructor{LNeg}} \coqdocvariable{A}) \coqdocvariable{l}\coqexternalref{:type scope:'exists' x '..' x ',' x}{http://coq.inria.fr/distrib/8.4pl3/stdlib/Coq.Init.Logic}{\coqdocnotation{)}}.\coqdoceol
\coqdocnoindent
\coqdockw{Definition} \coqdef{c completeness.complete mod.Valid CNF}{Valid\_CNF}{\coqdocdefinition{Valid\_CNF}} \coqdocvar{ll} := \coqdockw{\ensuremath{\forall}} \coqdocvar{l}, \coqexternalref{In}{http://coq.inria.fr/distrib/8.4pl3/stdlib/Coq.Lists.List}{\coqdocdefinition{In}} \coqdocvariable{l} \coqdocvariable{ll}\ensuremath{\rightarrow}\coqref{c completeness.complete mod.Valid Clause}{\coqdocdefinition{Valid\_Clause}} \coqdocvariable{l}.\coqdoceol
\end{code}
\subsection{Formalized Proof}\label{sec_complete_proof}
We now turn to prove the 6 statements described in Section \ref{sec_complete_sketch}.
\begin{code}
\coqdocnoindent
\coqdockw{Lemma} \coqdef{c completeness.complete mod.NNF equiv valid}{NNF\_equiv\_valid}{\coqdoclemma{NNF\_equiv\_valid}} : \coqdockw{\ensuremath{\forall}} \coqdocvar{v} \coqdocvar{A}, \coqref{c completeness.Y.TrueQ}{\coqdocdefinition{TrueQ}} \coqdocvariable{v} (\coqref{c completeness.complete mod.NNFtoPropF}{\coqdocdefinition{NNFtoPropF}} (\coqref{c completeness.complete mod.MakeNNF}{\coqdocdefinition{MakeNNF}}  \coqdocvariable{A}))\coqexternalref{:type scope:x '=' x}{http://coq.inria.fr/distrib/8.4pl3/stdlib/Coq.Init.Logic}{\coqdocnotation{=}}\coqref{c completeness.Y.TrueQ}{\coqdocdefinition{TrueQ}} \coqdocvariable{v}  \coqdocvariable{A} \coqexternalref{:type scope:x '/x5C' x}{http://coq.inria.fr/distrib/8.4pl3/stdlib/Coq.Init.Logic}{\coqdocnotation{\ensuremath{\land}}}\coqdoceol
\coqdocnoindent\hphantom{\coqdockw{Lemma} \coqdef{c completeness.complete mod.NNF equiv valid}{NNF\_equiv\_valid}{\coqdoclemma{NNF\_equiv\_valid}} : \coqdockw{\ensuremath{\forall}} \coqdocvar{v} \coqdocvar{A}, }%
\coqref{c completeness.Y.TrueQ}{\coqdocdefinition{TrueQ}} \coqdocvariable{v} (\coqref{c completeness.complete mod.NNFtoPropF}{\coqdocdefinition{NNFtoPropF}} (\coqref{c completeness.complete mod.MakeNNFN}{\coqdocdefinition{MakeNNFN}} \coqdocvariable{A}))\coqexternalref{:type scope:x '=' x}{http://coq.inria.fr/distrib/8.4pl3/stdlib/Coq.Init.Logic}{\coqdocnotation{=}}\coqref{c completeness.Y.TrueQ}{\coqdocdefinition{TrueQ}} \coqdocvariable{v} \coqref{b soundness.sound mod.:My scope:'xC2xAC' x}{\coqdocnotation{¬}}\coqdocvariable{A}.\coqdoceol
\coqdocnoindent
\coqdockw{Theorem} \coqdef{c completeness.complete mod.CNF equiv valid}{CNF\_equiv\_valid}{\coqdoclemma{CNF\_equiv\_valid}} : \coqdockw{\ensuremath{\forall}} \coqdocvar{v} \coqdocvar{A}, \coqref{c completeness.Y.TrueQ}{\coqdocdefinition{TrueQ}} \coqdocvariable{v} (\coqref{c completeness.complete mod.CNFtoPropF}{\coqdocdefinition{CNFtoPropF}} (\coqref{c completeness.complete mod.MakeCNF}{\coqdocdefinition{MakeCNF}} \coqdocvariable{A})) \coqexternalref{:type scope:x '=' x}{http://coq.inria.fr/distrib/8.4pl3/stdlib/Coq.Init.Logic}{\coqdocnotation{=}} \coqref{c completeness.Y.TrueQ}{\coqdocdefinition{TrueQ}} \coqdocvariable{v} (\coqref{c completeness.complete mod.NNFtoPropF}{\coqdocdefinition{NNFtoPropF}} \coqdocvariable{A}).\coqdoceol
\coqdocnoindent
\coqdockw{Theorem} \coqdef{c completeness.complete mod.CNF valid}{CNF\_valid}{\coqdoclemma{CNF\_valid}} : \coqdockw{\ensuremath{\forall}} \coqdocvar{ll}, \coqref{c completeness.Y.Valid}{\coqdocdefinition{Valid}} (\coqref{c completeness.complete mod.CNFtoPropF}{\coqdocdefinition{CNFtoPropF}} \coqdocvariable{ll}) \ensuremath{\rightarrow} \coqref{c completeness.complete mod.Valid CNF}{\coqdocdefinition{Valid\_CNF}} \coqdocvariable{ll}.\coqdoceol
\coqdocnoindent
\coqdockw{Theorem} \coqdef{c completeness.complete mod.CNF provable}{CNF\_provable}{\coqdoclemma{CNF\_provable}} : \coqdockw{\ensuremath{\forall}} \coqdocvar{ll}, \coqref{c completeness.complete mod.Valid CNF}{\coqdocdefinition{Valid\_CNF}} \coqdocvariable{ll} \ensuremath{\rightarrow} \coqref{c completeness.Y.Provable}{\coqdocdefinition{Provable}} (\coqref{c completeness.complete mod.CNFtoPropF}{\coqdocdefinition{CNFtoPropF}} \coqdocvariable{ll}).\coqdoceol
\coqdocnoindent
\coqdockw{Theorem} \coqdef{c completeness.complete mod.CNF impl prov}{CNF\_impl\_prov}{\coqdoclemma{CNF\_impl\_prov}} : \coqdockw{\ensuremath{\forall}} \coqdocvar{A}, \coqref{c completeness.Y.Provable}{\coqdocdefinition{Provable}} (\coqref{c completeness.complete mod.CNFtoPropF}{\coqdocdefinition{CNFtoPropF}} (\coqref{c completeness.complete mod.MakeCNF}{\coqdocdefinition{MakeCNF}} \coqdocvariable{A}) \coqref{b soundness.sound mod.:My scope:x 'xE2x86x92' x}{\coqdocnotation{→}} \coqref{c completeness.complete mod.NNFtoPropF}{\coqdocdefinition{NNFtoPropF}} \coqdocvariable{A}).\coqdoceol
\coqdocnoindent
\coqdockw{Lemma} \coqdef{c completeness.complete mod.NNF impl prov}{NNF\_impl\_prov}{\coqdoclemma{NNF\_impl\_prov}} : \coqdockw{\ensuremath{\forall}} \coqdocvar{A}, \coqref{c completeness.Y.Provable}{\coqdocdefinition{Provable}} (\coqref{c completeness.complete mod.NNFtoPropF}{\coqdocdefinition{NNFtoPropF}} (\coqref{c completeness.complete mod.MakeNNF}{\coqdocdefinition{MakeNNF}}  \coqdocvariable{A}) \coqref{b soundness.sound mod.:My scope:x 'xE2x86x92' x}{\coqdocnotation{→}}  \coqdocvariable{A}) \coqexternalref{:type scope:x '/x5C' x}{http://coq.inria.fr/distrib/8.4pl3/stdlib/Coq.Init.Logic}{\coqdocnotation{\ensuremath{\land}}}\coqdoceol
\coqdocnoindent\hphantom{\coqdockw{Lemma} \coqdef{c completeness.complete mod.NNF impl prov}{NNF\_impl\_prov}{\coqdoclemma{NNF\_impl\_prov}} : \coqdockw{\ensuremath{\forall}} \coqdocvar{A}, }%
\coqref{c completeness.Y.Provable}{\coqdocdefinition{Provable}} (\coqref{c completeness.complete mod.NNFtoPropF}{\coqdocdefinition{NNFtoPropF}} (\coqref{c completeness.complete mod.MakeNNFN}{\coqdocdefinition{MakeNNFN}} \coqdocvariable{A}) \coqref{b soundness.sound mod.:My scope:x 'xE2x86x92' x}{\coqdocnotation{→}} \coqref{b soundness.sound mod.:My scope:'xC2xAC' x}{\coqdocnotation{¬}}\coqdocvariable{A}).\coqdoceol
\end{code}
In a previous version I proved \coqdoclemma{CNF\_valid} with contraposition as described above. In that version I needed to first prove that syntactical validity was decidable. In the latest version I prove this directly, by showing that a clause is valid if it is valid under the countervaluation for it (which assigns true to a variable $p$ iff $\neg p$ occurs in the clause).

This allows us to prove completeness:
\begin{code}
\coqdocnoindent
\coqdockw{Theorem} \coqdef{c completeness.complete mod.Completeness}{Completeness}{\coqdoclemma{Completeness}} : \coqref{c completeness.Y.Prop Completeness}{\coqdocdefinition{Prop\_Completeness}}.\coqdoceol
\end{code}

\section{Equivalence between the calculi}\label{sec_equiv}
We now turn to formalizing the equivalence between the natural deduction calculus Nc, Hilbert-style calculus Hc and Gentzen's sequent calculus Gc.
\subsection{Nc and Hc}\label{sec_equiv_Hc}
In \verb"d_hilbert_calculus" I define Hilbert-style calculus. Natural deduction calculus and Hilbert-style calculus are quite similar. The difference is that in Hilbert-style calculus we have axioms, which allows us to only have implication elimination (modus ponens) as inference rule. The Hilbert-style calculus is defined below.
\begin{code}
\coqdocnoindent
\coqdockw{Inductive} \coqdef{d hilbert calculus.hilbert mod.AxiomH}{AxiomH}{\coqdocinductive{AxiomH}} : \coqref{d hilbert calculus.S.PropF}{\coqdocinductive{PropF}} \ensuremath{\rightarrow} \coqdockw{Prop} :=\coqdoceol
\coqdocnoindent
\ensuremath{|} \coqdef{d hilbert calculus.hilbert mod.HOrI1}{HOrI1}{\coqdocconstructor{HOrI1}}  : \coqdockw{\ensuremath{\forall}} \coqdocvar{A} \coqdocvar{B}  , \coqref{d hilbert calculus.AxiomH}{\coqdocinductive{AxiomH}} (\coqdocvariable{A} \coqref{b soundness.sound mod.:My scope:x 'xE2x86x92' x}{\coqdocnotation{→}} \coqdocvariable{A}\coqref{b soundness.sound mod.:My scope:x 'xE2x88xA8' x}{\coqdocnotation{∨}}\coqdocvariable{B})\coqdoceol
\coqdocnoindent
\ensuremath{|} \coqdef{d hilbert calculus.hilbert mod.HOrI2}{HOrI2}{\coqdocconstructor{HOrI2}}  : \coqdockw{\ensuremath{\forall}} \coqdocvar{A} \coqdocvar{B}  , \coqref{d hilbert calculus.AxiomH}{\coqdocinductive{AxiomH}} (\coqdocvariable{B} \coqref{b soundness.sound mod.:My scope:x 'xE2x86x92' x}{\coqdocnotation{→}} \coqdocvariable{A}\coqref{b soundness.sound mod.:My scope:x 'xE2x88xA8' x}{\coqdocnotation{∨}}\coqdocvariable{B})\coqdoceol
\coqdocnoindent
\ensuremath{|} \coqdef{d hilbert calculus.hilbert mod.HAndI}{HAndI}{\coqdocconstructor{HAndI}}  : \coqdockw{\ensuremath{\forall}} \coqdocvar{A} \coqdocvar{B}  , \coqref{d hilbert calculus.AxiomH}{\coqdocinductive{AxiomH}} (\coqdocvariable{A} \coqref{b soundness.sound mod.:My scope:x 'xE2x86x92' x}{\coqdocnotation{→}} \coqdocvariable{B} \coqref{b soundness.sound mod.:My scope:x 'xE2x86x92' x}{\coqdocnotation{→}} \coqdocvariable{A}\coqref{b soundness.sound mod.:My scope:x 'xE2x88xA7' x}{\coqdocnotation{∧}}\coqdocvariable{B})\coqdoceol
\coqdocnoindent
\ensuremath{|} \coqdef{d hilbert calculus.hilbert mod.HOrE}{HOrE}{\coqdocconstructor{HOrE}}   : \coqdockw{\ensuremath{\forall}} \coqdocvar{A} \coqdocvar{B} \coqdocvar{C}, \coqref{d hilbert calculus.AxiomH}{\coqdocinductive{AxiomH}} (\coqdocvariable{A}\coqref{b soundness.sound mod.:My scope:x 'xE2x88xA8' x}{\coqdocnotation{∨}}\coqdocvariable{B} \coqref{b soundness.sound mod.:My scope:x 'xE2x86x92' x}{\coqdocnotation{→}} \coqref{b soundness.sound mod.:My scope:x 'xE2x86x92' x}{\coqdocnotation{(}}\coqdocvariable{A} \coqref{b soundness.sound mod.:My scope:x 'xE2x86x92' x}{\coqdocnotation{→}} \coqdocvariable{C}\coqref{b soundness.sound mod.:My scope:x 'xE2x86x92' x}{\coqdocnotation{)}} \coqref{b soundness.sound mod.:My scope:x 'xE2x86x92' x}{\coqdocnotation{→}} \coqref{b soundness.sound mod.:My scope:x 'xE2x86x92' x}{\coqdocnotation{(}}\coqdocvariable{B} \coqref{b soundness.sound mod.:My scope:x 'xE2x86x92' x}{\coqdocnotation{→}} \coqdocvariable{C}\coqref{b soundness.sound mod.:My scope:x 'xE2x86x92' x}{\coqdocnotation{)}} \coqref{b soundness.sound mod.:My scope:x 'xE2x86x92' x}{\coqdocnotation{→}} \coqdocvariable{C})\coqdoceol
\coqdocnoindent
\ensuremath{|} \coqdef{d hilbert calculus.hilbert mod.HAndE1}{HAndE1}{\coqdocconstructor{HAndE1}} : \coqdockw{\ensuremath{\forall}} \coqdocvar{A} \coqdocvar{B}  , \coqref{d hilbert calculus.AxiomH}{\coqdocinductive{AxiomH}} (\coqdocvariable{A}\coqref{b soundness.sound mod.:My scope:x 'xE2x88xA7' x}{\coqdocnotation{∧}}\coqdocvariable{B} \coqref{b soundness.sound mod.:My scope:x 'xE2x86x92' x}{\coqdocnotation{→}} \coqdocvariable{A})\coqdoceol
\coqdocnoindent
\ensuremath{|} \coqdef{d hilbert calculus.hilbert mod.HAndE2}{HAndE2}{\coqdocconstructor{HAndE2}} : \coqdockw{\ensuremath{\forall}} \coqdocvar{A} \coqdocvar{B}  , \coqref{d hilbert calculus.AxiomH}{\coqdocinductive{AxiomH}} (\coqdocvariable{A}\coqref{b soundness.sound mod.:My scope:x 'xE2x88xA7' x}{\coqdocnotation{∧}}\coqdocvariable{B} \coqref{b soundness.sound mod.:My scope:x 'xE2x86x92' x}{\coqdocnotation{→}} \coqdocvariable{B})\coqdoceol
\coqdocnoindent
\ensuremath{|} \coqdef{d hilbert calculus.hilbert mod.HS}{HS}{\coqdocconstructor{HS}}      : \coqdockw{\ensuremath{\forall}} \coqdocvar{A} \coqdocvar{B} \coqdocvar{C}, \coqref{d hilbert calculus.AxiomH}{\coqdocinductive{AxiomH}} (\coqref{b soundness.sound mod.:My scope:x 'xE2x86x92' x}{\coqdocnotation{(}}\coqdocvariable{A} \coqref{b soundness.sound mod.:My scope:x 'xE2x86x92' x}{\coqdocnotation{→}} \coqdocvariable{B} \coqref{b soundness.sound mod.:My scope:x 'xE2x86x92' x}{\coqdocnotation{→}} \coqdocvariable{C}\coqref{b soundness.sound mod.:My scope:x 'xE2x86x92' x}{\coqdocnotation{)}} \coqref{b soundness.sound mod.:My scope:x 'xE2x86x92' x}{\coqdocnotation{→}} \coqref{b soundness.sound mod.:My scope:x 'xE2x86x92' x}{\coqdocnotation{(}}\coqdocvariable{A} \coqref{b soundness.sound mod.:My scope:x 'xE2x86x92' x}{\coqdocnotation{→}} \coqdocvariable{B}\coqref{b soundness.sound mod.:My scope:x 'xE2x86x92' x}{\coqdocnotation{)}} \coqref{b soundness.sound mod.:My scope:x 'xE2x86x92' x}{\coqdocnotation{→}} \coqdocvariable{A} \coqref{b soundness.sound mod.:My scope:x 'xE2x86x92' x}{\coqdocnotation{→}} \coqdocvariable{C})\coqdoceol
\coqdocnoindent
\ensuremath{|} \coqdef{d hilbert calculus.hilbert mod.HK}{HK}{\coqdocconstructor{HK}}      : \coqdockw{\ensuremath{\forall}} \coqdocvar{A} \coqdocvar{B}  , \coqref{d hilbert calculus.AxiomH}{\coqdocinductive{AxiomH}} (\coqdocvariable{A} \coqref{b soundness.sound mod.:My scope:x 'xE2x86x92' x}{\coqdocnotation{→}} \coqdocvariable{B} \coqref{b soundness.sound mod.:My scope:x 'xE2x86x92' x}{\coqdocnotation{→}} \coqdocvariable{A})\coqdoceol
\coqdocnoindent
\ensuremath{|} \coqdef{d hilbert calculus.hilbert mod.HClas}{HClas}{\coqdocconstructor{HClas}}  : \coqdockw{\ensuremath{\forall}} \coqdocvar{A}    , \coqref{d hilbert calculus.AxiomH}{\coqdocinductive{AxiomH}} (\coqref{b soundness.sound mod.:My scope:'xC2xAC' x}{\coqdocnotation{¬(¬}}\coqdocvariable{A}\coqref{b soundness.sound mod.:My scope:'xC2xAC' x}{\coqdocnotation{)}} \coqref{b soundness.sound mod.:My scope:x 'xE2x86x92' x}{\coqdocnotation{→}} \coqdocvariable{A})\coqdoceol
\coqdocnoindent
.\coqdoceol
\coqdocnoindent
\coqdockw{Inductive} \coqdef{d hilbert calculus.hilbert mod.Hc}{Hc}{\coqdocinductive{Hc}} : \coqexternalref{list}{http://coq.inria.fr/distrib/8.4pl3/stdlib/Coq.Init.Datatypes}{\coqdocinductive{list}} \coqref{d hilbert calculus.S.PropF}{\coqdocinductive{PropF}}\ensuremath{\rightarrow} \coqref{d hilbert calculus.S.PropF}{\coqdocinductive{PropF}}\ensuremath{\rightarrow}\coqdockw{Prop} :=\coqdoceol
\coqdocnoindent
\ensuremath{|} \coqdef{d hilbert calculus.hilbert mod.Hass}{Hass}{\coqdocconstructor{Hass}}  : \coqdockw{\ensuremath{\forall}} \coqdocvar{A} \coqdocvar{Γ},   \coqexternalref{In}{http://coq.inria.fr/distrib/8.4pl3/stdlib/Coq.Lists.List}{\coqdocdefinition{In}} \coqdocvariable{A} \coqdocvariable{Γ}               \ensuremath{\rightarrow} \coqdocvariable{Γ} \coqref{d hilbert calculus.hilbert mod.:My scope:x 'xE2x8AxA2H' x}{\coqdocnotation{⊢}}\coqref{d hilbert calculus.hilbert mod.:My scope:x 'xE2x8AxA2H' x}{\coqdocnotation{H}} \coqdocvariable{A}\coqdoceol
\coqdocnoindent
\ensuremath{|} \coqdef{d hilbert calculus.hilbert mod.Hax}{Hax}{\coqdocconstructor{Hax}}   : \coqdockw{\ensuremath{\forall}} \coqdocvar{A} \coqdocvar{Γ},   \coqref{d hilbert calculus.hilbert mod.AxiomH}{\coqdocinductive{AxiomH}} \coqdocvariable{A}             \ensuremath{\rightarrow} \coqdocvariable{Γ} \coqref{d hilbert calculus.hilbert mod.:My scope:x 'xE2x8AxA2H' x}{\coqdocnotation{⊢}}\coqref{d hilbert calculus.hilbert mod.:My scope:x 'xE2x8AxA2H' x}{\coqdocnotation{H}} \coqdocvariable{A}\coqdoceol
\coqdocnoindent
\ensuremath{|} \coqdef{d hilbert calculus.hilbert mod.HImpE}{HImpE}{\coqdocconstructor{HImpE}} : \coqdockw{\ensuremath{\forall}} \coqdocvar{Γ} \coqdocvar{A} \coqdocvar{B}, \coqdocvariable{Γ} \coqref{d hilbert calculus.hilbert mod.:My scope:x 'xE2x8AxA2H' x}{\coqdocnotation{⊢}}\coqref{d hilbert calculus.hilbert mod.:My scope:x 'xE2x8AxA2H' x}{\coqdocnotation{H}} \coqdocvariable{A} \coqref{b soundness.sound mod.:My scope:x 'xE2x86x92' x}{\coqdocnotation{→}} \coqdocvariable{B} \ensuremath{\rightarrow} \coqdocvariable{Γ} \coqref{d hilbert calculus.hilbert mod.:My scope:x 'xE2x8AxA2H' x}{\coqdocnotation{⊢}}\coqref{d hilbert calculus.hilbert mod.:My scope:x 'xE2x8AxA2H' x}{\coqdocnotation{H}} \coqdocvariable{A} \ensuremath{\rightarrow} \coqdocvariable{Γ} \coqref{d hilbert calculus.hilbert mod.:My scope:x 'xE2x8AxA2H' x}{\coqdocnotation{⊢}}\coqref{d hilbert calculus.hilbert mod.:My scope:x 'xE2x8AxA2H' x}{\coqdocnotation{H}} \coqdocvariable{B}\coqdoceol
\coqdocnoindent
\coqdockw{where} \coqdef{d hilbert calculus.hilbert mod.:My scope:x 'xE2x8AxA2H' x}{"}{"}Γ ⊢H A" := (\coqref{d hilbert calculus.Hc}{\coqdocinductive{Hc}} \coqdocvar{Γ} \coqdocvar{A}) : \coqdocvar{My\_scope}.\coqdoceol
\end{code}
I then prove the equivalence of these systems.
\begin{code}
\coqdocnoindent
\coqdockw{Theorem} \coqdef{d hilbert calculus.hilbert mod.Nc equiv Hc}{Nc\_equiv\_Hc}{\coqdoclemma{Nc\_equiv\_Hc}} : \coqdockw{\ensuremath{\forall}} \coqdocvar{Γ} \coqdocvar{A}, \coqdocvariable{Γ} \coqref{b soundness.sound mod.:My scope:x 'xE2x8AxA2' x}{\coqdocnotation{⊢}} \coqdocvariable{A} \coqexternalref{:type scope:x '<->' x}{http://coq.inria.fr/distrib/8.4pl3/stdlib/Coq.Init.Logic}{\coqdocnotation{\ensuremath{\leftrightarrow}}} \coqdocvariable{Γ} \coqref{d hilbert calculus.hilbert mod.:My scope:x 'xE2x8AxA2H' x}{\coqdocnotation{⊢}}\coqref{d hilbert calculus.hilbert mod.:My scope:x 'xE2x8AxA2H' x}{\coqdocnotation{H}} \coqdocvariable{A}.\coqdoceol
\end{code}
The proof is not hard. In the direction from left to right, one basically needs to prove that all axioms for Hc are provable in Nc. In the converse direction, one needs to prove the Deduction Theorem for Hc, which states that
$$\Gamma,A\vdash_{\Hc} B\iff\Gamma\vdash_{\Hc} A\to B.$$
\subsection{Nc and Gc}\label{sec_equiv_Gc}
In \verb"e_sequent_calculus" I define Gentzen's sequent calculus. Gentzen's sequent calculus is quite different than natural deduction. In sequent calculus the propositions are \emph{sequents} of the form $\Gamma\vdash\Delta$, where $\Gamma$ and $\Delta$ are either sequences, multisets or sets of formulae. A sequent intuitively means that the conjunction of $\Gamma$ implies the disjunction of $\Delta$. Instead of having introduction and elimination rules as in natural deduction, there are left rules and right rules for each connective, which introduce that connective on that side of the sequent.

I define the sequent calculus in the formalization below. Since lists are easiest to work with, I use them in sequents. Usually when lists are used in sequent calculus, there are some structural rules allowing to move formulas around in the list. However, having these rules make induction proofs very hard, and requires more nonstructural induction proofs, so instead I defined the sequent calculus in such a way that it is possible to apply a rule anywhere in the lists (ordinarily the relevant formulas appear only on one end of the lists). I use the symbol $\supset$ as separator in a sequent.
\begin{code}
\coqdocnoindent
\coqdockw{Inductive} \coqdef{e sequent calculus.sequent mod.G}{G}{\coqdocinductive{G}} : \coqexternalref{list}{http://coq.inria.fr/distrib/8.4pl3/stdlib/Coq.Init.Datatypes}{\coqdocinductive{list}} \coqref{e sequent calculus.S.PropF}{\coqdocinductive{PropF}}\ensuremath{\rightarrow}\coqexternalref{list}{http://coq.inria.fr/distrib/8.4pl3/stdlib/Coq.Init.Datatypes}{\coqdocinductive{list}} \coqref{e sequent calculus.S.PropF}{\coqdocinductive{PropF}}\ensuremath{\rightarrow}\coqdockw{Prop} :=\coqdoceol
\coqdocnoindent
\ensuremath{|} \coqdef{e sequent calculus.sequent mod.Gax}{Gax}{\coqdocconstructor{Gax}}  : \coqdockw{\ensuremath{\forall}} \coqdocvar{A} \coqdocvar{Γ} \coqdocvar{Δ}      , \coqexternalref{In}{http://coq.inria.fr/distrib/8.4pl3/stdlib/Coq.Lists.List}{\coqdocdefinition{In}} \coqdocvariable{A} \coqdocvariable{Γ}           \ensuremath{\rightarrow} \coqexternalref{In}{http://coq.inria.fr/distrib/8.4pl3/stdlib/Coq.Lists.List}{\coqdocdefinition{In}} \coqdocvariable{A} \coqdocvariable{Δ}          \ensuremath{\rightarrow} \coqdocvariable{Γ} \coqref{e sequent calculus.sequent mod.:My scope:x 'xE2x8Ax83' x}{\coqdocnotation{⊃}} \coqdocvariable{Δ}\coqdoceol
\coqdocnoindent
\ensuremath{|} \coqdef{e sequent calculus.sequent mod.GBot}{GBot}{\coqdocconstructor{GBot}} : \coqdockw{\ensuremath{\forall}} \coqdocvar{Γ} \coqdocvar{Δ}        , \coqexternalref{In}{http://coq.inria.fr/distrib/8.4pl3/stdlib/Coq.Lists.List}{\coqdocdefinition{In}} \coqref{b soundness.sound mod.:My scope:'xE2x8AxA5'}{\coqdocnotation{⊥}} \coqdocvariable{Γ}                              \ensuremath{\rightarrow} \coqdocvariable{Γ} \coqref{e sequent calculus.sequent mod.:My scope:x 'xE2x8Ax83' x}{\coqdocnotation{⊃}} \coqdocvariable{Δ}\coqdoceol
\coqdocnoindent
\ensuremath{|} \coqdef{e sequent calculus.sequent mod.AndL}{AndL}{\coqdocconstructor{AndL}} : \coqdockw{\ensuremath{\forall}} \coqdocvar{A} \coqdocvar{B} \coqdocvar{Γ1} \coqdocvar{Γ2} \coqdocvar{Δ}, \coqdocvariable{Γ1}\coqexternalref{:list scope:x '++' x}{http://coq.inria.fr/distrib/8.4pl3/stdlib/Coq.Init.Datatypes}{\coqdocnotation{++}}\coqdocvariable{A}\coqexternalref{:list scope:x '::' x}{http://coq.inria.fr/distrib/8.4pl3/stdlib/Coq.Init.Datatypes}{\coqdocnotation{::}}\coqdocvariable{B}\coqexternalref{:list scope:x '::' x}{http://coq.inria.fr/distrib/8.4pl3/stdlib/Coq.Init.Datatypes}{\coqdocnotation{::}}\coqdocvariable{Γ2} \coqref{e sequent calculus.sequent mod.:My scope:x 'xE2x8Ax83' x}{\coqdocnotation{⊃}} \coqdocvariable{Δ}                    \ensuremath{\rightarrow} \coqdocvariable{Γ1}\coqexternalref{:list scope:x '++' x}{http://coq.inria.fr/distrib/8.4pl3/stdlib/Coq.Init.Datatypes}{\coqdocnotation{++}}\coqdocvariable{A}\coqref{b soundness.sound mod.:My scope:x 'xE2x88xA7' x}{\coqdocnotation{∧}}\coqdocvariable{B}\coqexternalref{:list scope:x '::' x}{http://coq.inria.fr/distrib/8.4pl3/stdlib/Coq.Init.Datatypes}{\coqdocnotation{::}}\coqdocvariable{Γ2} \coqref{e sequent calculus.sequent mod.:My scope:x 'xE2x8Ax83' x}{\coqdocnotation{⊃}} \coqdocvariable{Δ}\coqdoceol
\coqdocnoindent
\ensuremath{|} \coqdef{e sequent calculus.sequent mod.AndR}{AndR}{\coqdocconstructor{AndR}} : \coqdockw{\ensuremath{\forall}} \coqdocvar{A} \coqdocvar{B} \coqdocvar{Γ} \coqdocvar{Δ1} \coqdocvar{Δ2}, \coqdocvariable{Γ} \coqref{e sequent calculus.sequent mod.:My scope:x 'xE2x8Ax83' x}{\coqdocnotation{⊃}} \coqdocvariable{Δ1}\coqexternalref{:list scope:x '++' x}{http://coq.inria.fr/distrib/8.4pl3/stdlib/Coq.Init.Datatypes}{\coqdocnotation{++}}\coqdocvariable{A}\coqexternalref{:list scope:x '::' x}{http://coq.inria.fr/distrib/8.4pl3/stdlib/Coq.Init.Datatypes}{\coqdocnotation{::}}\coqdocvariable{Δ2}    \ensuremath{\rightarrow} \coqdocvariable{Γ} \coqref{e sequent calculus.sequent mod.:My scope:x 'xE2x8Ax83' x}{\coqdocnotation{⊃}} \coqdocvariable{Δ1}\coqexternalref{:list scope:x '++' x}{http://coq.inria.fr/distrib/8.4pl3/stdlib/Coq.Init.Datatypes}{\coqdocnotation{++}}\coqdocvariable{B}\coqexternalref{:list scope:x '::' x}{http://coq.inria.fr/distrib/8.4pl3/stdlib/Coq.Init.Datatypes}{\coqdocnotation{::}}\coqdocvariable{Δ2}   \ensuremath{\rightarrow} \coqdocvariable{Γ} \coqref{e sequent calculus.sequent mod.:My scope:x 'xE2x8Ax83' x}{\coqdocnotation{⊃}} \coqdocvariable{Δ1}\coqexternalref{:list scope:x '++' x}{http://coq.inria.fr/distrib/8.4pl3/stdlib/Coq.Init.Datatypes}{\coqdocnotation{++}}\coqdocvariable{A}\coqref{b soundness.sound mod.:My scope:x 'xE2x88xA7' x}{\coqdocnotation{∧}}\coqdocvariable{B}\coqexternalref{:list scope:x '::' x}{http://coq.inria.fr/distrib/8.4pl3/stdlib/Coq.Init.Datatypes}{\coqdocnotation{::}}\coqdocvariable{Δ2}\coqdoceol
\coqdocnoindent
\ensuremath{|} \coqdef{e sequent calculus.sequent mod.OrL}{OrL}{\coqdocconstructor{OrL}}  : \coqdockw{\ensuremath{\forall}} \coqdocvar{A} \coqdocvar{B} \coqdocvar{Γ1} \coqdocvar{Γ2} \coqdocvar{Δ}, \coqdocvariable{Γ1}\coqexternalref{:list scope:x '++' x}{http://coq.inria.fr/distrib/8.4pl3/stdlib/Coq.Init.Datatypes}{\coqdocnotation{++}}\coqdocvariable{A}\coqexternalref{:list scope:x '::' x}{http://coq.inria.fr/distrib/8.4pl3/stdlib/Coq.Init.Datatypes}{\coqdocnotation{::}}\coqdocvariable{Γ2} \coqref{e sequent calculus.sequent mod.:My scope:x 'xE2x8Ax83' x}{\coqdocnotation{⊃}} \coqdocvariable{Δ}    \ensuremath{\rightarrow} \coqdocvariable{Γ1}\coqexternalref{:list scope:x '++' x}{http://coq.inria.fr/distrib/8.4pl3/stdlib/Coq.Init.Datatypes}{\coqdocnotation{++}}\coqdocvariable{B}\coqexternalref{:list scope:x '::' x}{http://coq.inria.fr/distrib/8.4pl3/stdlib/Coq.Init.Datatypes}{\coqdocnotation{::}}\coqdocvariable{Γ2} \coqref{e sequent calculus.sequent mod.:My scope:x 'xE2x8Ax83' x}{\coqdocnotation{⊃}} \coqdocvariable{Δ}   \ensuremath{\rightarrow} \coqdocvariable{Γ1}\coqexternalref{:list scope:x '++' x}{http://coq.inria.fr/distrib/8.4pl3/stdlib/Coq.Init.Datatypes}{\coqdocnotation{++}}\coqdocvariable{A}\coqref{b soundness.sound mod.:My scope:x 'xE2x88xA8' x}{\coqdocnotation{∨}}\coqdocvariable{B}\coqexternalref{:list scope:x '::' x}{http://coq.inria.fr/distrib/8.4pl3/stdlib/Coq.Init.Datatypes}{\coqdocnotation{::}}\coqdocvariable{Γ2} \coqref{e sequent calculus.sequent mod.:My scope:x 'xE2x8Ax83' x}{\coqdocnotation{⊃}} \coqdocvariable{Δ}\coqdoceol
\coqdocnoindent
\ensuremath{|} \coqdef{e sequent calculus.sequent mod.OrR}{OrR}{\coqdocconstructor{OrR}}  : \coqdockw{\ensuremath{\forall}} \coqdocvar{A} \coqdocvar{B} \coqdocvar{Γ} \coqdocvar{Δ1} \coqdocvar{Δ2}, \coqdocvariable{Γ} \coqref{e sequent calculus.sequent mod.:My scope:x 'xE2x8Ax83' x}{\coqdocnotation{⊃}} \coqdocvariable{Δ1}\coqexternalref{:list scope:x '++' x}{http://coq.inria.fr/distrib/8.4pl3/stdlib/Coq.Init.Datatypes}{\coqdocnotation{++}}\coqdocvariable{A}\coqexternalref{:list scope:x '::' x}{http://coq.inria.fr/distrib/8.4pl3/stdlib/Coq.Init.Datatypes}{\coqdocnotation{::}}\coqdocvariable{B}\coqexternalref{:list scope:x '::' x}{http://coq.inria.fr/distrib/8.4pl3/stdlib/Coq.Init.Datatypes}{\coqdocnotation{::}}\coqdocvariable{Δ2}                     \ensuremath{\rightarrow} \coqdocvariable{Γ} \coqref{e sequent calculus.sequent mod.:My scope:x 'xE2x8Ax83' x}{\coqdocnotation{⊃}} \coqdocvariable{Δ1}\coqexternalref{:list scope:x '++' x}{http://coq.inria.fr/distrib/8.4pl3/stdlib/Coq.Init.Datatypes}{\coqdocnotation{++}}\coqdocvariable{A}\coqref{b soundness.sound mod.:My scope:x 'xE2x88xA8' x}{\coqdocnotation{∨}}\coqdocvariable{B}\coqexternalref{:list scope:x '::' x}{http://coq.inria.fr/distrib/8.4pl3/stdlib/Coq.Init.Datatypes}{\coqdocnotation{::}}\coqdocvariable{Δ2}\coqdoceol
\coqdocnoindent
\ensuremath{|} \coqdef{e sequent calculus.sequent mod.ImpL}{ImpL}{\coqdocconstructor{ImpL}} : \coqdockw{\ensuremath{\forall}} \coqdocvar{A} \coqdocvar{B} \coqdocvar{Γ1} \coqdocvar{Γ2} \coqdocvar{Δ}, \coqdocvariable{Γ1}\coqexternalref{:list scope:x '++' x}{http://coq.inria.fr/distrib/8.4pl3/stdlib/Coq.Init.Datatypes}{\coqdocnotation{++}}\coqdocvariable{B}\coqexternalref{:list scope:x '::' x}{http://coq.inria.fr/distrib/8.4pl3/stdlib/Coq.Init.Datatypes}{\coqdocnotation{::}}\coqdocvariable{Γ2} \coqref{e sequent calculus.sequent mod.:My scope:x 'xE2x8Ax83' x}{\coqdocnotation{⊃}} \coqdocvariable{Δ}    \ensuremath{\rightarrow} \coqdocvariable{Γ1}\coqexternalref{:list scope:x '++' x}{http://coq.inria.fr/distrib/8.4pl3/stdlib/Coq.Init.Datatypes}{\coqdocnotation{++}}\coqdocvariable{Γ2} \coqref{e sequent calculus.sequent mod.:My scope:x 'xE2x8Ax83' x}{\coqdocnotation{⊃}} \coqdocvariable{A}\coqexternalref{:list scope:x '::' x}{http://coq.inria.fr/distrib/8.4pl3/stdlib/Coq.Init.Datatypes}{\coqdocnotation{::}}\coqdocvariable{Δ}   \ensuremath{\rightarrow} \coqdocvariable{Γ1}\coqexternalref{:list scope:x '++' x}{http://coq.inria.fr/distrib/8.4pl3/stdlib/Coq.Init.Datatypes}{\coqdocnotation{++}}\coqdocvariable{A}\coqref{b soundness.sound mod.:My scope:x 'xE2x86x92' x}{\coqdocnotation{→}}\coqdocvariable{B}\coqexternalref{:list scope:x '::' x}{http://coq.inria.fr/distrib/8.4pl3/stdlib/Coq.Init.Datatypes}{\coqdocnotation{::}}\coqdocvariable{Γ2} \coqref{e sequent calculus.sequent mod.:My scope:x 'xE2x8Ax83' x}{\coqdocnotation{⊃}} \coqdocvariable{Δ}\coqdoceol
\coqdocnoindent
\ensuremath{|} \coqdef{e sequent calculus.sequent mod.ImpR}{ImpR}{\coqdocconstructor{ImpR}} : \coqdockw{\ensuremath{\forall}} \coqdocvar{A} \coqdocvar{B} \coqdocvar{Γ} \coqdocvar{Δ1} \coqdocvar{Δ2}, \coqdocvariable{A}\coqexternalref{:list scope:x '::' x}{http://coq.inria.fr/distrib/8.4pl3/stdlib/Coq.Init.Datatypes}{\coqdocnotation{::}}\coqdocvariable{Γ} \coqref{e sequent calculus.sequent mod.:My scope:x 'xE2x8Ax83' x}{\coqdocnotation{⊃}} \coqdocvariable{Δ1}\coqexternalref{:list scope:x '++' x}{http://coq.inria.fr/distrib/8.4pl3/stdlib/Coq.Init.Datatypes}{\coqdocnotation{++}}\coqdocvariable{B}\coqexternalref{:list scope:x '::' x}{http://coq.inria.fr/distrib/8.4pl3/stdlib/Coq.Init.Datatypes}{\coqdocnotation{::}}\coqdocvariable{Δ2}                     \ensuremath{\rightarrow} \coqdocvariable{Γ} \coqref{e sequent calculus.sequent mod.:My scope:x 'xE2x8Ax83' x}{\coqdocnotation{⊃}} \coqdocvariable{Δ1}\coqexternalref{:list scope:x '++' x}{http://coq.inria.fr/distrib/8.4pl3/stdlib/Coq.Init.Datatypes}{\coqdocnotation{++}}\coqdocvariable{A}\coqref{b soundness.sound mod.:My scope:x 'xE2x86x92' x}{\coqdocnotation{→}}\coqdocvariable{B}\coqexternalref{:list scope:x '::' x}{http://coq.inria.fr/distrib/8.4pl3/stdlib/Coq.Init.Datatypes}{\coqdocnotation{::}}\coqdocvariable{Δ2}\coqdoceol
\coqdocnoindent
\ensuremath{|} \coqdef{e sequent calculus.sequent mod.Cut}{Cut}{\coqdocconstructor{Cut}}  : \coqdockw{\ensuremath{\forall}} \coqdocvar{A} \coqdocvar{Γ} \coqdocvar{Δ}      , \coqdocvariable{Γ} \coqref{e sequent calculus.sequent mod.:My scope:x 'xE2x8Ax83' x}{\coqdocnotation{⊃}} \coqdocvariable{A}\coqexternalref{:list scope:x '::' x}{http://coq.inria.fr/distrib/8.4pl3/stdlib/Coq.Init.Datatypes}{\coqdocnotation{::}}\coqdocvariable{Δ}         \ensuremath{\rightarrow} \coqdocvariable{A}\coqexternalref{:list scope:x '::' x}{http://coq.inria.fr/distrib/8.4pl3/stdlib/Coq.Init.Datatypes}{\coqdocnotation{::}}\coqdocvariable{Γ} \coqref{e sequent calculus.sequent mod.:My scope:x 'xE2x8Ax83' x}{\coqdocnotation{⊃}} \coqdocvariable{Δ}        \ensuremath{\rightarrow} \coqdocvariable{Γ} \coqref{e sequent calculus.sequent mod.:My scope:x 'xE2x8Ax83' x}{\coqdocnotation{⊃}} \coqdocvariable{Δ}\coqdoceol
\coqdocnoindent
\coqdockw{where} \coqdef{e sequent calculus.sequent mod.:My scope:x 'xE2x8Ax83' x}{"}{"}Γ ⊃ Δ" := (\coqref{e sequent calculus.G}{\coqdocinductive{G}} \coqdocvar{Γ} \coqdocvar{Δ}) : \coqdocvar{My\_scope}.\coqdoceol
\end{code}
We now turn to proving the equivalence between Nc and Gc. Neither direction is easy.

In the direction from Gc to Nc we prove that
$$\Gamma\supset\Delta\implies\Gamma\vdash\bigvee\Delta,$$
where $\bigvee\Delta$ is the disjunction of the formulae in $\Delta$ (similar to a clause, except that the entries can now be any formulae instead of only literals). However, proving this directly is very hard, since in the definition of sequents we make changes in the middle of the disjunction. This means that in every step we have to use multiple or-eliminations and then or-introductions, which gets really ugly. Instead, we prove the following intermediate lemma:
$$\Gamma\supset\Delta\implies\Gamma,\neg\Delta\vdash\bot,$$
where $\neg\Delta$ is the element-wise negation of formulae in $\Delta$. In the formalization this looks as follows.

\begin{code}
\coqdocnoindent
\coqdockw{Definition} \coqdef{e sequent calculus.sequent mod.BigOr}{BigOr}{\coqdocdefinition{BigOr}} := \coqexternalref{fold right}{http://coq.inria.fr/distrib/8.4pl3/stdlib/Coq.Lists.List}{\coqdocdefinition{fold\_right}} \coqref{e sequent calculus.S.Disj}{\coqdocconstructor{Disj}} \coqref{b soundness.sound mod.:My scope:'xE2x8AxA5'}{\coqdocnotation{⊥}}.\coqdoceol
\coqdocnoindent
\coqdockw{Notation} \coqdef{e sequent calculus.sequent mod.::'xE2x8Bx81' x}{"}{"}⋁ Δ" := (\coqref{e sequent calculus.sequent mod.BigOr}{\coqdocdefinition{BigOr}} \coqdocvar{Δ}) (\coqdoctac{at} \coqdockw{level} 19).\coqdoceol
\coqdockw{Notation} \coqdef{e sequent calculus.sequent mod.::'xC2xACl' x}{"}{"}¬l Γ" := (\coqexternalref{map}{http://coq.inria.fr/distrib/8.4pl3/stdlib/Coq.Lists.List}{\coqdocdefinition{map}} \coqref{e sequent calculus.S.Neg}{\coqdocdefinition{Neg}} \coqdocvar{Γ}) (\coqdoctac{at} \coqdockw{level} 40).\coqdoceol
\coqdocnoindent
\coqdockw{Lemma} \coqdef{e sequent calculus.sequent mod.G to Nc Neg}{G\_to\_Nc\_Neg}{\coqdoclemma{G\_to\_Nc\_Neg}} : \coqdockw{\ensuremath{\forall}} \coqdocvar{Γ} \coqdocvar{Δ}, \coqdocvariable{Γ} \coqref{e sequent calculus.sequent mod.:My scope:x 'xE2x8Ax83' x}{\coqdocnotation{⊃}} \coqdocvariable{Δ} \ensuremath{\rightarrow} \coqdocvariable{Γ}\coqexternalref{:list scope:x '++' x}{http://coq.inria.fr/distrib/8.4pl3/stdlib/Coq.Init.Datatypes}{\coqdocnotation{++}}\coqref{e sequent calculus.sequent mod.::'xC2xACl' x}{\coqdocnotation{¬}}\coqref{e sequent calculus.sequent mod.::'xC2xACl' x}{\coqdocnotation{l}} \coqdocvariable{Δ} \coqref{b soundness.sound mod.:My scope:x 'xE2x8AxA2' x}{\coqdocnotation{⊢}} \coqref{b soundness.sound mod.:My scope:'xE2x8AxA5'}{\coqdocnotation{⊥}}.\coqdoceol
\coqdockw{Theorem} \coqdef{e sequent calculus.sequent mod.G to Nc}{G\_to\_Nc}{\coqdoclemma{G\_to\_Nc}} : \coqdockw{\ensuremath{\forall}} \coqdocvar{Γ} \coqdocvar{Δ}, \coqdocvariable{Γ} \coqref{e sequent calculus.sequent mod.:My scope:x 'xE2x8Ax83' x}{\coqdocnotation{⊃}} \coqdocvariable{Δ} \ensuremath{\rightarrow} \coqdocvariable{Γ} \coqref{b soundness.sound mod.:My scope:x 'xE2x8AxA2' x}{\coqdocnotation{⊢}} \coqref{e sequent calculus.sequent mod.::'xE2x8Bx81' x}{\coqdocnotation{⋁}}\coqdocvariable{Δ}.\coqdoceol
\end{code}
In the direction from Gc to Nc the major difficulty is proving weakening for Gc. On paper proving weakening is not very hard, but in the formalization the fact that we can apply rules anywhere in the list, and the fact that the weakening also occurs anywhere in the list, makes it quite hard to prove this (for example, we need to distinguish cases in every induction step which of the two formulae occurs earlier in the list).
\begin{code}
\coqdocnoindent
\coqdockw{Lemma} \coqdef{e sequent calculus.sequent mod.WeakL}{WeakL}{\coqdoclemma{WeakL}} : \coqdockw{\ensuremath{\forall}} \coqdocvar{Γ1} \coqdocvar{Γ2} \coqdocvar{Δ} \coqdocvar{A}, \coqdocvariable{Γ1}\coqexternalref{:list scope:x '++' x}{http://coq.inria.fr/distrib/8.4pl3/stdlib/Coq.Init.Datatypes}{\coqdocnotation{++}}\coqdocvariable{Γ2} \coqref{e sequent calculus.sequent mod.:My scope:x 'xE2x8Ax83' x}{\coqdocnotation{⊃}} \coqdocvariable{Δ} \ensuremath{\rightarrow} \coqdocvariable{Γ1}\coqexternalref{:list scope:x '++' x}{http://coq.inria.fr/distrib/8.4pl3/stdlib/Coq.Init.Datatypes}{\coqdocnotation{++}}\coqdocvariable{A}\coqexternalref{:list scope:x '::' x}{http://coq.inria.fr/distrib/8.4pl3/stdlib/Coq.Init.Datatypes}{\coqdocnotation{::}}\coqdocvariable{Γ2} \coqref{e sequent calculus.sequent mod.:My scope:x 'xE2x8Ax83' x}{\coqdocnotation{⊃}} \coqdocvariable{Δ}.\coqdoceol
\coqdocnoindent
\coqdockw{Lemma} \coqdef{e sequent calculus.sequent mod.WeakR}{WeakR}{\coqdoclemma{WeakR}} : \coqdockw{\ensuremath{\forall}} \coqdocvar{Γ} \coqdocvar{Δ1} \coqdocvar{Δ2} \coqdocvar{A}, \coqdocvariable{Γ} \coqref{e sequent calculus.sequent mod.:My scope:x 'xE2x8Ax83' x}{\coqdocnotation{⊃}} \coqdocvariable{Δ1}\coqexternalref{:list scope:x '++' x}{http://coq.inria.fr/distrib/8.4pl3/stdlib/Coq.Init.Datatypes}{\coqdocnotation{++}}\coqdocvariable{Δ2} \ensuremath{\rightarrow} \coqdocvariable{Γ} \coqref{e sequent calculus.sequent mod.:My scope:x 'xE2x8Ax83' x}{\coqdocnotation{⊃}} \coqdocvariable{Δ1}\coqexternalref{:list scope:x '++' x}{http://coq.inria.fr/distrib/8.4pl3/stdlib/Coq.Init.Datatypes}{\coqdocnotation{++}}\coqdocvariable{A}\coqexternalref{:list scope:x '::' x}{http://coq.inria.fr/distrib/8.4pl3/stdlib/Coq.Init.Datatypes}{\coqdocnotation{::}}\coqdocvariable{Δ2}.\coqdoceol
\coqdocnoindent
\coqdockw{Theorem} \coqdef{e sequent calculus.sequent mod.Nc to G}{Nc\_to\_G}{\coqdoclemma{Nc\_to\_G}} : \coqdockw{\ensuremath{\forall}} \coqdocvar{Γ} \coqdocvar{A}, \coqdocvariable{Γ} \coqref{b soundness.sound mod.:My scope:x 'xE2x8AxA2' x}{\coqdocnotation{⊢}} \coqdocvariable{A} \ensuremath{\rightarrow} \coqdocvariable{Γ} \coqref{e sequent calculus.sequent mod.:My scope:x 'xE2x8Ax83' x}{\coqdocnotation{⊃}} \coqexternalref{ListNotations.:list scope:'[' x ';' '..' ';' x ']'}{http://coq.inria.fr/distrib/8.4pl3/stdlib/Coq.Lists.List}{\coqdocnotation{[}}\coqdocvariable{A}\coqexternalref{ListNotations.:list scope:'[' x ';' '..' ';' x ']'}{http://coq.inria.fr/distrib/8.4pl3/stdlib/Coq.Lists.List}{\coqdocnotation{]}}.\coqdoceol
\coqdocemptyline
\coqdocnoindent
\end{code}

\section{Cut Elimination}\label{sec_cut_elim}
In \verb"f_cut_elimination" I prove the cut elimination theorem. I first introduce the cut-free sequent calculus, which has no cut rule. To strengthen the cut elimination a bit, I only allow the axiom rule for atomic formulae in the cut-free calculus. I denote the cut-free sequents as $\Gamma\,\coqdocnotation{⊃c}\,\Delta$ (I don't give the definition here; it is very similar to the definition \coqdocinductive{G}).

I give a semantic proof for cut elimination. This proof does not generalize easily to predicate calculus, the semantic proof for predicate calculus is way more involved. I also tried to give a syntactical proof, but failed due to time contraints.

The semantic proof goes as follows. First we prove soundness of the sequent calculus with cut, which states that the conclusion of any derivation (with cuts) is valid. The main part is to prove completeness of the cut-free calculus. Suppose that we are given a valid sequent. If no logical connectives occur in the sequent, then we prove that we can either apply the axiom rule or the falsum rule, using that the sequent is valid. On the other hand, if any connective occurs in the sequent, we apply the corresponding rule (in reverse) to eliminate that connective. We can then show the hypotheses to the rule are also valid, and they contain fewer connectives than the conclusion. So we can use a nonstructural induction to conclude that the hypotheses are provable, which finishes the proof.

We need the following definitions. The first line defines the notation for a valid sequent. 
Then we define the size of a formula and sequent, which in this case is the number of connectives occurring in it.
\begin{code}
\coqdocnoindent
\coqdockw{Definition} \coqdef{c completeness.complete mod.Validates}{Validates}{\coqdocdefinition{Validates}} \coqdocvar{v} \coqdocvar{Δ} := \coqexternalref{:type scope:'exists' x '..' x ',' x}{http://coq.inria.fr/distrib/8.4pl3/stdlib/Coq.Init.Logic}{\coqdocnotation{\ensuremath{\exists}}} \coqdocvar{A}\coqexternalref{:type scope:'exists' x '..' x ',' x}{http://coq.inria.fr/distrib/8.4pl3/stdlib/Coq.Init.Logic}{\coqdocnotation{,}} \coqexternalref{In}{http://coq.inria.fr/distrib/8.4pl3/stdlib/Coq.Lists.List}{\coqdocdefinition{In}} \coqdocvariable{A} \coqdocvariable{Δ} \coqexternalref{:type scope:x '/x5C' x}{http://coq.inria.fr/distrib/8.4pl3/stdlib/Coq.Init.Logic}{\coqdocnotation{\ensuremath{\land}}} \coqexternalref{Is true}{http://coq.inria.fr/distrib/8.4pl3/stdlib/Coq.Bool.Bool}{\coqdocdefinition{Is\_true}} (\coqref{c completeness.Y.TrueQ}{\coqdocdefinition{TrueQ}} \coqdocvariable{v} \coqdocvariable{A}).\coqdoceol
\coqdocnoindent
\coqdockw{Notation} \coqdef{f cut elimination.cut mod.semantic proof.::x '=xE2x8Ax83' x}{"}{"}Γ =⊃ Δ" := (\coqdockw{\ensuremath{\forall}} \coqdocvar{v},\coqref{f cut elimination.S.Satisfies}{\coqdocdefinition{Satisfies}} \coqdocvariable{v} \coqdocvar{Γ}\ensuremath{\rightarrow}\coqref{f cut elimination.cut mod.Validates}{\coqdocdefinition{Validates}} \coqdocvariable{v} \coqdocvar{Δ}) (\coqdoctac{at} \coqdockw{level} 80).\coqdoceol
\coqdocnoindent
\coqdockw{Fixpoint} \coqdef{f cut elimination.cut mod.size}{size}{\coqdocdefinition{size}} \coqdocvar{A} : \coqexternalref{nat}{http://coq.inria.fr/distrib/8.4pl3/stdlib/Coq.Init.Datatypes}{\coqdocinductive{nat}} := \coqdockw{match} \coqdocvariable{A} \coqdockw{with}\coqdoceol
\coqdocindent{0.50em}
\ensuremath{|} \coqref{b soundness.sound mod.:My scope:'x23' x}{\coqdocnotation{\#}} \coqdocvar{P}   \ensuremath{\Rightarrow} 0\coqdoceol
\coqdocindent{0.50em}
\ensuremath{|} \coqref{b soundness.sound mod.:My scope:'xE2x8AxA5'}{\coqdocnotation{⊥}}     \ensuremath{\Rightarrow} 0\coqdoceol
\coqdocindent{0.50em}
\ensuremath{|} \coqdocvar{B} \coqref{b soundness.sound mod.:My scope:x 'xE2x88xA8' x}{\coqdocnotation{∨}} \coqdocvar{C} \ensuremath{\Rightarrow} \coqexternalref{S}{http://coq.inria.fr/distrib/8.4pl3/stdlib/Coq.Init.Datatypes}{\coqdocconstructor{S}} (\coqref{f cut elimination.size}{\coqdocdefinition{size}} \coqdocvar{B} \coqexternalref{:nat scope:x '+' x}{http://coq.inria.fr/distrib/8.4pl3/stdlib/Coq.Init.Peano}{\coqdocnotation{+}} \coqref{f cut elimination.size}{\coqdocdefinition{size}} \coqdocvar{C})\coqdoceol
\coqdocindent{0.50em}
\ensuremath{|} \coqdocvar{B} \coqref{b soundness.sound mod.:My scope:x 'xE2x88xA7' x}{\coqdocnotation{∧}} \coqdocvar{C} \ensuremath{\Rightarrow} \coqexternalref{S}{http://coq.inria.fr/distrib/8.4pl3/stdlib/Coq.Init.Datatypes}{\coqdocconstructor{S}} (\coqref{f cut elimination.size}{\coqdocdefinition{size}} \coqdocvar{B} \coqexternalref{:nat scope:x '+' x}{http://coq.inria.fr/distrib/8.4pl3/stdlib/Coq.Init.Peano}{\coqdocnotation{+}} \coqref{f cut elimination.size}{\coqdocdefinition{size}} \coqdocvar{C})\coqdoceol
\coqdocindent{0.50em}
\ensuremath{|} \coqdocvar{B} \coqref{b soundness.sound mod.:My scope:x 'xE2x86x92' x}{\coqdocnotation{→}} \coqdocvar{C} \ensuremath{\Rightarrow} \coqexternalref{S}{http://coq.inria.fr/distrib/8.4pl3/stdlib/Coq.Init.Datatypes}{\coqdocconstructor{S}} (\coqref{f cut elimination.size}{\coqdocdefinition{size}} \coqdocvar{B} \coqexternalref{:nat scope:x '+' x}{http://coq.inria.fr/distrib/8.4pl3/stdlib/Coq.Init.Peano}{\coqdocnotation{+}} \coqref{f cut elimination.size}{\coqdocdefinition{size}} \coqdocvar{C})\coqdoceol
\coqdocnoindent
\coqdockw{end}.\coqdoceol
\coqdocnoindent
\coqdockw{Definition} \coqdef{f cut elimination.cut mod.sizel}{sizel}{\coqdocdefinition{sizel}} := \coqref{f cut elimination.B.map fold right}{\coqdocdefinition{map\_fold\_right}} \coqref{f cut elimination.cut mod.size}{\coqdocdefinition{size}} \coqexternalref{plus}{http://coq.inria.fr/distrib/8.4pl3/stdlib/Coq.Init.Peano}{\coqdocdefinition{plus}} 0.\coqdoceol
\coqdocnoindent
\coqdockw{Definition} \coqdef{f cut elimination.cut mod.sizes}{sizes}{\coqdocdefinition{sizes}} \coqdocvar{Γ} \coqdocvar{Δ}:= \coqref{f cut elimination.cut mod.sizel}{\coqdocdefinition{sizel}} \coqdocvariable{Γ} \coqexternalref{:nat scope:x '+' x}{http://coq.inria.fr/distrib/8.4pl3/stdlib/Coq.Init.Peano}{\coqdocnotation{+}} \coqref{f cut elimination.cut mod.sizel}{\coqdocdefinition{sizel}} \coqdocvariable{Δ}.\coqdoceol
\end{code}
This allows us to prove the following theorem by induction to $n$.
\begin{code}
\coqdocnoindent
\coqdockw{Theorem} \coqdef{f cut elimination.cut mod.Gcf complete induction}{Gcf\_complete\_induction}{\coqdoclemma{Gcf\_complete\_induction}} : \coqdockw{\ensuremath{\forall}} \coqdocvar{n} \coqdocvar{Γ} \coqdocvar{Δ}, \coqref{f cut elimination.cut mod.sizes}{\coqdocdefinition{sizes}} \coqdocvariable{Γ} \coqdocvariable{Δ} \coqexternalref{:nat scope:x '<=' x}{http://coq.inria.fr/distrib/8.4pl3/stdlib/Coq.Init.Peano}{\coqdocnotation{\ensuremath{\le}}} \coqdocvariable{n} \ensuremath{\rightarrow} \coqdocvariable{Γ} \coqref{f cut elimination.cut mod.semantic proof.::x '=xE2x8Ax83' x}{\coqdocnotation{=⊃}} \coqdocvariable{Δ} \ensuremath{\rightarrow} \coqdocvariable{Γ} \coqref{f cut elimination.cut mod.semantic proof.:My scope:x 'xE2x8Ax83c' x}{\coqdocnotation{⊃}}\coqref{f cut elimination.cut mod.semantic proof.:My scope:x 'xE2x8Ax83c' x}{\coqdocnotation{c}} \coqdocvariable{Δ}.\coqdoceol
\end{code}
Together with Soundness this proves cut elimination.
\begin{code}
\coqdocnoindent
\coqdockw{Theorem} \coqdef{f cut elimination.cut mod.Cut elimination}{Cut\_elimination}{\coqdoclemma{Cut\_elimination}} : \coqdockw{\ensuremath{\forall}} \coqdocvar{Γ} \coqdocvar{Δ}, \coqdocvariable{Γ} \coqref{e sequent calculus.sequent mod.:My scope:x 'xE2x8Ax83' x}{\coqdocnotation{⊃}} \coqdocvariable{Δ} \ensuremath{\rightarrow} \coqdocvariable{Γ} \coqref{f cut elimination.cut mod.semantic proof.:My scope:x 'xE2x8Ax83c' x}{\coqdocnotation{⊃}}\coqref{f cut elimination.cut mod.semantic proof.:My scope:x 'xE2x8Ax83c' x}{\coqdocnotation{c}} \coqdocvariable{Δ}.\coqdoceol
\end{code}
\section{Conclusion and Future work}
In this project I have shown that the proofs of basic but important theorems about propositional calculus can be implemented in a proof assistant relatively painlessly. Since many concepts in proof theory are naturally definable using induction, Coq's excellent support for inductive datatypes and proofs by induction makes the formalization of this subject relatively easy. The lack of strong automation in Coq does require that the user has to focus a lot on details in the proof.

This project can be seen as a proof of concept for formalizing theorems in proof theory. One can extend this project in numerous ways. One can consider nonclassical logics, and prove soundness and completeness results about them, like the completeness intuitionistic logic w.r.t. Kripke models. Another way to extend this is to consider predicate calculus instead of propositional calculus. This has been done in the literature already, for example to prove Gödel's first incompleteness theorem.\citep{OConnor2005Incompleteness} Another possible extension is to prove other theorems, like the Beth definability theorem or Craig's interpolation theorem.

\addcontentsline{toc}{section}{References}
\bibliography{references}

\end{document}